\documentclass[opre,nonblindrev]{informs4}

\OneAndAHalfSpacedXI 
\usepackage{natbib}
 \bibpunct[, ]{(}{)}{,}{a}{}{,}%
 \def\bibfont{\small}%
 \def\bibsep{\smallskipamount}%
\TheoremsNumberedThrough    
\ECRepeatTheorems
\EquationsNumberedThrough    
\MANUSCRIPTNO{}

\usepackage{algorithm} 
\usepackage{algorithmic} 
\usepackage{tikz} 
\usepackage{enumerate}
\usepackage{rotating,lineno}
\usepackage{flafter}
\usepackage[T1]{fontenc}
\usepackage[utf8]{inputenc}
\usepackage{babel}
\usepackage[font=small,labelfont=bf]{caption}
\usepackage{graphicx}
\usepackage{capt-of}
\usepackage[mathscr]{euscript}
\usepackage{pgfplots}
\usepackage{multirow}
\usepackage{array, makecell}
\usepackage{hyperref}
\hypersetup{
    colorlinks=true,
    linkcolor=blue,
    filecolor=magenta,      
    urlcolor=blue,
    citecolor={blue},
}
\pgfplotsset{compat=1.12}

\DeclareSymbolFont{rsfs}{U}{rsfs}{m}{n}
\DeclareSymbolFontAlphabet{\mathscrsfs}{rsfs}

\definecolor{lightgray}{RGB}{229, 231, 233}

\usepackage{amsmath}

\usepackage{nicefrac}
\usepackage[shortlabels]{enumitem}

\newcommand{\bA}{ \mathbf{A} }
\newcommand{\ba}{ \mathbf{a} }
\newcommand{\bb}{ \mathbf{b} }
\newcommand{\bB}{ \mathbf{B} }

\newcommand{\bc}{ \mathbf{c} }
\newcommand{\bD}{ \mathbf{D} }
\newcommand{\bd}{ \mathbf{d} }

\newcommand{\bE}{ \mathbf{E} }
\newcommand{\besmall}{ \mathbf{e} }

\newcommand{\bG}{ {\mathbf{A}}_{\cT} }

\newcommand{\bh}{ {\mathbf{b}}_{\cT} }

\newcommand{\bP}{ \mathbf{P} }

\newcommand{\bS}{ \mathbf{S} }

\newcommand{\bu}{ \Bar{\mathbf{y}} }
\newcommand{\bv}{ \mathbf{v} }

\newcommand{\bX}{ \mathbf{X} }
\newcommand{\bx}{ \mathbf{x} }

\newcommand{\by}{ \mathbf{y} }

\newcommand{\bz}{ \mathbf{z} }

\newcommand{\bone}{\mathbf{1}}
\newcommand{\bzero}{\mathbf{0}}

\newcommand{\FO}{\text{FO}}
\newcommand{\IO}{\text{IL}} 
\newcommand{\CIO}{\text{IO}}

\newcommand{\IL}{\text{IL}}

\newcommand{\MIIO}{\text{M2}}
\newcommand{\GIL}{\text{IL$^{\emptyset}$}}
\newcommand{\MGIL}{\text{GIL}}
\newcommand{\GGIL}{\text{GIL$_r$}}
\newcommand{\ADG}{\text{M1}}
\newcommand{\bepsilon}{\mathbf{\epsilon}}

\newcommand{\desirable}{relevant}
\newcommand{\Desirable}{Relevant}
\newcommand{\objective}{trivial}
\newcommand{\Objective}{Trivial}

\newcommand{\pref}{preferred}
\newcommand{\Pref}{Preferred}

\newcommand{\DFO}{\text{DFO}}

\newcommand{\cJ}{\mathcal{J}}
\newcommand{\cK}{\mathcal{K}}

\newcommand{\cT}{\mathcal{T}}
\newcommand{\cR}{\mathcal{R}}
\newcommand{\cP}{\mathcal{P}}

\newcommand{\cF}{ \mathscr{F} }
\newcommand{\cD}{ \mathscr{D} }
\newcommand{\cQ}{ \mathscr{Q} }
\newcommand{\ObservationConstraintTradeoff}{observation-goal tradeoff}

\tikzstyle{very loosely dotted}=          [dash pattern=on \pgflinewidth off 6pt]
\tikzstyle{very loosely dashdotted}=      [dash pattern=on 3pt off 8pt on \the\pgflinewidth off 8pt]
\tikzstyle{costumdashed}=                  [dash pattern=on 10pt off 3pt]

\newcommand{\fa}[1]{\textit{\textcolor{black!40!green}{#1}}}

\newcommand{\todo}[1]{\textit{\textcolor{red}{#1}}}

\begin{document}

\RUNTITLE{Inverse Learning}

\TITLE{Inverse Learning: Solving Partially Known Models Using Inverse Optimization}

\ARTICLEAUTHORS{%
\AUTHOR{Farzin Ahmadi, Fardin Ganjkhanloo, Kimia Ghobadi}
\AFF{Department of Civil and Systems Engineering, The Center for Systems Science and Engineering, The Malone Center for Engineering in Healthcare, Johns Hopkins University, Baltimore, MD 21218, \EMAIL{fahmadi1@jhu.edu}, \EMAIL{fganjkh1@jhu.edu}, \EMAIL{kimia@jhu.edu}} 

} 
\RUNAUTHOR{Ahmadi, Ganjkhanloo, Ghobadi}
\ABSTRACT{

We consider the problem of learning optimal solutions of a partially known linear optimization problem and recovering its underlying cost function where a set of past decisions and the feasible set are known. We develop a new framework, denoted as \emph{Inverse Learning}, that extends the inverse optimization literature to (1) learn the optimal solution of the underlying problem, (2) integrate additional information on constraints and their importance, and (3) control the balance between mimicking past behaviors and reaching new goals and rules for the learned solution. We pose inverse learning as an optimization problem that maps given (feasible and infeasible) observations to a single optimal solution with minimum perturbation, hence, not only recovering the missing cost vector but also providing an optimal solution simultaneously. 
The framework provides insights into an essential tradeoff in recovering linear optimization problems with regard to preserving observed behaviors and binding constraints of the known feasible set at optimality. We propose a series of mixed integer linear programming models to capture the effects of this tradeoff and validate it using a two-dimensional example. We then demonstrate the framework's applicability to a diet recommendation problem for a population of hypertension and prediabetic patients. The goal is to balance dietary constraints to achieve the necessary nutritional goals with the dietary habits of the users to encourage adherence to the diet.  Results indicate that our models recommend daily food intakes that preserve the original data trends while providing a range of options to patients and providers based on the tradeoff. 

}

\maketitle

\section{Introduction}

Optimization models are used in solving many real-world problems, from healthcare applications in diet recommendations and cancer treatment planning to energy and service industries. However, deriving the necessary parameters to write these optimization models remains largely non-intuitive despite frequent access to historical observations of solutions to the problems.  
One method that is capable of recovering unknown parameters of partially known optimization problems is Inverse Optimization \citep{Ahuja01}. Inverse optimization recovers the unknown parameters of a partially known (forward) optimization problem based on given observed solution(s).  
Recent studies have shown that inverse optimization models can successfully recover missing parameters of forward linear optimization problems ranging from linear optimization problems to integer or convex optimization models with a focus on recovering the unknown cost function~\citep{beil2003inverse, iyengar2005inverse, bertsimas2012inverse, chan2014generalized, chan2019inverse,aswani2018inverse,shahmoradi2021quantile}. 

In this work, we develop a novel framework, denoted as \emph{Inverse Learning}, that recovers unknown cost functions of an underlying forward optimization problem based on given observation(s), similar to most inverse optimization models. 
The framework is, however, structurally distinct from inverse optimization paradigms in three key features of 
(1) directly learning an optimal solution to the underlying forward optimization problem in addition to recovering the missing cost function, hence, enabling preservation of the characteristics of observed solutions,  
(2) integrating constraint priority knowledge (if any) into the inverse framework to  further guide the optimal solutions' recovery, 
and (3) providing a range of solutions along with a control mechanism to balance the inherent tradeoff between optimal solutions that retain observations' characteristics (feature 1) and those that bind prioritized constraints and goals (feature 2). 

The first key feature of the inverse learning framework is unifying the process of learning forward optimal solutions with the recovery of the underlying cost vector. This process is in contrast with most inverse optimization paradigms, where the primary goal is to recover a cost vector to the forward problem. Then, they often solve the recovered forward optimization problem to find an optimal solution \citep{Ahuja01, keshavarz2011imputing, chan2019inverse}. 
Decoupling of finding the forward cost vectors and optimal solutions leads to several disadvantages, including diminished control over the characteristics of the obtained optimal solution, instability in the optimal solutions and cost vectors, a higher computational cost of the inverse model, and the need to solve the forward model independently. 
Simultaneous recovery of both the optimal solutions and the cost vectors enables the inverse framework to better control the characteristics of eventual optimal solutions and cost vectors instead of indirectly impacting the solutions through cost vectors. The additional knowledge of solutions' properties can be embedded into the models, along with any other available information about the cost vector. 
The inverse literature considers various methods to choose the most appropriate cost vector, 
including secondary objectives, prior beliefs, or side constraints \citep{ghobadi2018robust, chan2014generalized,shahmoradi2021quantile}. 
However, these methods only apply to cost vectors and do not guarantee that the desired characteristics will translate to optimal solutions, which may be of more interest to users. 
The emphasis on primarily recovering cost vectors may also lead to instability of the obtained optimal solutions. 
For example, the Simplex method returns extreme point(s) as optimal solutions that may be far from the observed solutions, as highlighted in detail by \cite{shahmoradi2021quantile}. 
Additionally, concurrent recovery of the optimal solutions and the cost vectors reduces computational complexity by reducing the number of non-convex constraints in the inverse model and removing the need to solve the forward problem independently.  


Affording additional flexibility in the model can integrate data-driven or user-provided knowledge in the inverse model. 
Figure \ref{fig:ILvsSIO_comparison} shows a schematic example for three observations (gray dots) and a feasible region (shaded area). The recovered cost vector and optimal solution are shown for inverse optimization (left) and inverse learning (right). A sequential recovery of first the cost vector and then the optimal solution, often through Simplex-based algorithms, may result in solutions that are far from the given observation, and little control is afforded on their characteristics once the cost vector is fixed. It may also lead to instability as either extreme point may be returned by the forward model. Simultaneous recovery of both enables enhanced control over the solution properties and integration of desired characteristics, as the figure illustrates.

\begin{figure}[h]
\begin{center}
\begin{tikzpicture}[scale=0.5][
  bluenode/.style={shape=rectangle, draw=cyan, line width=2, fill= cyan},
  greennode/.style={shape=bar, draw=gray, line width=2, fill= gray},
  graynode/.style={shape=circle, draw=gray, line width=2, fill= gray},
  blacknode/.style={shape=circle, draw=black, line width=2, fill= black},
  40greennode/.style={shape=circle, draw=black!40!green, line width=2, fill= black!40!green},
  rednode/.style={shape=rectangle, draw=red, line width=2, fill= red}
]

\fill  [black!10!white] (5,1) -- (7,3) -- (7,5)-- (5,7) -- (3,7) -- (1,5) --cycle;
\draw [line width=1pt][gray][costumdashed] (5,1) -- (7,3);
\draw [line width=1pt][gray][costumdashed] (7,3) -- (7,5);
\draw [line width=1pt][gray][costumdashed] (7,5) -- (5,7);
\draw [line width=1pt][gray][costumdashed] (5,7) -- (2.95,7);
\draw [line width=1pt][gray][costumdashed] (3,7) -- (1,5);
\draw [line width=1pt][gray][costumdashed] (1,5) -- (5,1);

\draw [very thick][dashed][red] (4.8, 5.1) -- (5.9,6.1);
\draw[fill][gray] (4.8, 5.1) circle [radius=0.1];
\draw [very thick][dashed][red] (5.4, 5.1) -- (6.2,5.9);
\draw[fill][gray] (5.4, 5.1) circle [radius=0.1];
\draw [very thick][dashed][red] (6, 5.2) -- (6.5,5.7);
\draw[fill][gray] (6, 5.2) circle [radius=0.1];

\draw[fill][black!40!red] (5,7) circle [radius=0.15];
\draw [very thick][black!40!green] (6,7) -- (7,6);
\draw[very thick][->][black!40!green] (6.5,6.5)  -- (7,7) node[right]{};

\node[label={}: {Inverse Optimization}] at (5,-0.5)  {};
\end{tikzpicture}
\qquad
\begin{tikzpicture}[scale=0.5][
  bluenode/.style={shape=rectangle, draw=cyan, line width=2, fill= cyan},
  greennode/.style={shape=bar, draw=gray, line width=2, fill= gray},
  graynode/.style={shape=circle, draw=gray, line width=2, fill= gray},
  blacknode/.style={shape=circle, draw=black, line width=2, fill= black},
  greennode/.style={shape=circle, draw=black!40!green, line width=2, fill= black!40!green},
  rednode/.style={shape=rectangle, draw=red, line width=2, fill= red}
]

\fill  [black!10!white] (5,1) -- (7,3) -- (7,5)-- (5,7) -- (3,7) -- (1,5) --cycle;
\draw [line width=1][gray][costumdashed] (5,1) -- (7,3);
\draw [line width=1][gray][costumdashed] (7,3) -- (7,5);
\draw [line width=1][gray][costumdashed] (7,5) -- (5,7);
\draw [line width=1][gray][costumdashed] (5,7) -- (2.95,7);
\draw [line width=1][gray][costumdashed] (3,7) -- (1,5);
\draw [line width=1][gray][costumdashed] (1,5) -- (5,1);

\draw [very thick][dashed][red] (4.8, 5.1) -- (6.2,5.8);
\draw[fill][gray] (4.8, 5.1) circle [radius=0.1];
\draw [very thick][dashed][red] (5.4, 5.1) -- (6.2,5.8);
\draw[fill][gray] (5.4, 5.1) circle [radius=0.1];
\draw [very thick][dashed][red] (6, 5.2) -- (6.2,5.8);
\draw[fill][gray] (6, 5.2) circle [radius=0.1];
\draw[fill][black!40!green] (6.2,5.8) circle [radius=0.15];
\draw [very thick][black!40!green] (6,7) -- (7,6);
\draw[very thick][->][black!40!green] (6.5,6.5)  -- (7,7) node[right]{};

\node[label={}: {Inverse Learning}] at (5,-0.5)  {};
\end{tikzpicture}
\begin{tikzpicture}[scale=0.3][] 
\footnotesize
\draw [line width=4pt][gray] (9.4,15) -- (10,15); \node[label={right,label distance=0.1cm}: {Observations}] at (10.5,15)  {};
\draw [line width=4pt][black!40!green] (9.4,13.6) -- (10,13.6); \node[label={right,label distance=0.1cm}: {Inverse Output}] at (10.5,13.6)  {};
\draw [line width=4pt][black!40!red] (9.4,12.2) -- (10,12.2); \node[label={right,label distance=0.1cm}: {Forward Output}] at (10.5,12.2)  {};
\draw [line width=1pt][gray][costumdashed] (8.4,10.8) -- (11.4,10.8); \node[label={right,label distance=0.1cm}: {Constraint}] at (10.5,10.8)  {};
\draw [line width=4pt][white] (8.4,4) -- (9,4);
\end{tikzpicture}
\caption{Inverse optimization (left) and inverse learning (right) approaches are applied to a set of three observations (blue dots). Inverse optimization obtains the cost function (green) and often solves the recovered forward problem to find an optimal solution (red dot), which may or may not represent the observations. Inverse learning simultaneously detects the closest solution and a cost vector that renders it optimal, hence, consolidating additional observations' properties. 
}
\label{fig:ILvsSIO_comparison}
\end{center}
\end{figure}

The second feature of inverse learning is embedding additional constraints' priority knowledge in the learning process, should such information exists. 
Integration of observations into inverse optimization frameworks has been studied extensively in the literature, especially in the form of a distance function {\citep{chan2021inversesurvey,chan2014generalized, keshavarz2011imputing, aswani2018inverse, ghobadi2018robust}}. However, little attention has been paid to the information conveyed by constraints beyond identifying the feasible region. In particular, those constraints that are active at optimality play vital roles in determining the missing cost vectors and the eventual optimal solution of the forward problem. They may also carry significance in practical settings. For instance, in resource allocation problems, meeting all demands or utilizing the minimum amount of certain resources may be desirable, or in diet problems, maximizing protein or minimizing saturated fat may be meaningful goals for the dieter. This additional information and goals can be provided by end-users, domain experts, or derived from data. The inverse learning framework can incorporate such goals by binding preferred constraints and considering different tiers of constraints in the models. If no additional information exists, the inverse learning model operates as conventional models and treats all constraints equally.


The third feature of inverse learning lies in its ability to control tradeoffs between observation-driven (feature 1) and goal-driven (feature 2) outcomes and provide a decision-support system based on user's preferences or by progressively adding goals. The framework aims to balance between the two main inputs to inverse models: observations -- as a measure of past behavior or decisions -- and constraints -- as a measure of rules, restrictions, or goals. 
The framework can produce a range of recommendations that traverse between mimicking observations most closely 
to finding ideal solutions that bind the maximum number of desired constraints to resemble a goal-driven approach.  
The ability to provide a spectrum of solutions enables tailoring the solutions to the needs and preferences of the users. 
The span of solutions facilitates comparison between adhering to goals or maintaining resemblance to observations and provides a base to navigate the tradeoff between them. It can also be employed as a mechanism to carve a path for the user to gradually yet progressively reach additional goals by initially preserving past behavior (mimicking observations) and then successively augmenting more goals into the model. 

The focus on the learning process in inverse frameworks is motivated by practical settings where examples of historical solutions or observations are abundant but may not be optimal or even feasible for the underlying forward problem. For instance, expert decisions in healthcare planning or treatments, customer decisions in service industries, or patients' dietary decisions, among others. These solutions highlight the behavior and decisions of the user or the status of the system. In these settings, finding optimal solutions that represent observations may be paramount as they represent underlying patterns or preferences that may not be captured in the model otherwise. Incorporating inputs from observations and constraint priorities directly into the models guides finding optimal solutions that are relevant and appropriate for the problem. The additional flexibility and control also empower users to tailor the solutions or obtain a set of solutions that gradually improves goals, and indirectly, improve cost vector stability and sensitivity to outliers.

We demonstrate the applicability of the inverse learning framework on a diet recommendation problem for a population of 2,090 patients with a self-reported hypertension diagnosis. 
Precision nutrition is increasingly used as part of medical treatments to prevent or control symptoms in prediabetic, hypertension, or kidney failure patients. These patients are often prescribed a low-sodium diet as part of their treatment. A typical example of such a diet that clinicians recommend is the Dietary Approaches to Stop Hypertension (DASH) diet \citep{sacks2001effects,liese2009adherence} which aims to control sodium intake and is shown to lower blood pressure--a known risk factor of diabetes~\citep{sacks2001effects}. Studies suggest that adherence to the DASH diet has an inverse relation with type 2 diabetes and can reduce future risk of it by as much as 20\%  \citep{,liese2009adherence, Campbell2017DASHDiabetes}.  
However, the dietary guidelines provided to the patients are rarely personalized or tailored to individual patients' tastes or lifestyles, making long-term adherence to diets challenging \citep{downer2016predictors,inelmen2005predictors,bazrafkan2021overweight}. 
Observations of patients' past food intake can provide insights into their preferences, while the DASH diet constraints provide nutritional goals and constraints. A spectrum of diet recommendations that traverse from preserving patients' behavior to diets that increasingly improve nutritional goals (e.g., minimizing sodium, maximizing fiber, minimizing sugar, etc.) can provide patients with a range of choices. Starting with diets closer to patients' current behaviors and gradually introducing new goals may provide a more personalized diet recommendation that can improve patients' adherence as opposed to a one-size-fits-all approach.

We pose a forward optimization problem based on the DASH diet nutritional constraints with an unknown cost vector of patients' preferences and palate. We then use inverse learning models to recommend diets for patient classes. We demonstrate how inverse learning models can be used to directly recommend a set of diets with clear interpretations for the patients and the dietitians and recover the underlying cost function. The recommendations provide a control mechanism to balance the tradeoff between past food intake and nutritional goals. The added flexibility in recommending decisions advocates 
algorithm usability through algorithm aversion \citep{dietvorst2015algorithm,grand2022best} as it provides means for decision-makers and users to gain control over decisions resulting from the algorithms \citep{dietvorst2018overcoming}. The outcomes demonstrate a dietary path for the patients to initially maintain their habits as much as possible while adhering to the prescribed diet. Gradually, the recommended diets can bind more nutritional constraints and move towards meeting more dietary goals and the ideal nutrition plan for the patients. 

In the rest of this paper, we first provide an overview of the literature and then introduce the forward and inverse learning models in Section \ref{Section:Methodology}. We then extend to goal-integrated inverse learning models to control the \ObservationConstraintTradeoff in Section \ref{sec:AIL}. A detailed two-dimensional numerical example is then discussed in Section \ref{Section:NumericalEx}, and finally, we apply our full framework to the diet recommendation problem in Section \ref{Section:Application}. We showcase a web-based decision-support tool for this problem that enables exploring various diets based on the user's preferences and goals. Concluding remarks are provided in Section \ref{Section:Conclusion} and proofs are outlined in an electronic companion.

\subsection{Literature Review}


Early investigations of inverse problems can be traced back to the works of \cite{burton1992instance} who considered the inverse shortest path problem and discussed inverse combinatorial optimization problems. Inverse settings of combinatorial and general linear optimization problems gained more traction in  time \citep{yang1999two,zhang1996calculating,zhang1999further,zhang1999solution}.  \cite{heuberger2004inverse} provided a survey on inverse combinatorial problems and the existing methods for approaching such problems. However, \cite{Ahuja01} was the first to use the term ``Inverse Optimization'' and pose it as a general problem of recovering cost vectors of linear programming problems. Since then, a common setting in inverse optimization has been to recover the cost function and implicit preferences of a decision-maker in a constrained environment in the form of a (forward) optimization problem \citep{chan2014generalized,aswani2018inverse,esfahani2018data,shahmoradi2021quantile} and solve the recovered problem with known methods to find optimal solutions. While the majority of the literature considers the given observation(s) as a single optimal \citep{Ahuja01,iyengar2005inverse} or near-optimal point \citep{ghate2015inverse,chan2019inverse,zhang2010inverse, schaefer2009inverse, wang2009cutting,ghate2020inverse,lasserre2013inverse,naghavi2019inverse,roland2013inverse}, 
in practice, often a collection of past solutions can be observed from the forward problem \citep{keshavarz2011imputing,shahmoradi2021quantile,babier2021ensemble} or from similar but differently constrained problems \citep{vcerny2016inverse}. In such cases, existing methods rely on minimizing optimality gap measures in terms of primal-dual objective values to recover unknown parameters. However, it has been argued in the literature that such methods result in inconsistent estimations in data-driven settings \citep{aswani2018inverse}. 

Multi-observation settings in inverse optimization have increasingly gained traction due to their real-world applicability. Recent studies have considered various methods to recover linear problems in the presence of multiple observations for cost vectors \citep{babier2021ensemble, esfahani2018data, aswani2018inverse, tan2019deep, hewitt2020data, moghaddass2021inverse} or for constraint parameters \citep{ghobadi2021inferring,chan2020inverse,Schede2019learning}. Others have considered the cases in which a stream of observations is given, as opposed to an existing batch of observations \citep{dong2018generalized}. Statistical approaches in noisy settings \citep{aswani2019statistics, shahmoradi2021quantile,shahmoradi2022optimality} and deep learning models in data-driven settings are also explored \citep{tan2019deep}. 

Many applications of inverse optimization methodologies have also been studied including finance \citep{bertsimas2012inverse}, network and telecommunication \citep{farago2003inverse}, supply chain planning \citep{pibernik2011secure}, transportation problems \citep{chow2012inverse}, control systems \citep{9336679}, and auction mechanisms \citep{beil2003inverse}. 
\citet{barmann2018online} consider a setting where the decision-maker decides on the variables in different domains and provides an online-learning framework for inverse optimization. \cite{konstantakopoulos2017robust} combine inverse optimization and game theory methods to learn the utility of agents in a non-cooperative game. Other studies consider healthcare applications to employ inverse optimization's unique ability to recover unknown parameters in multi-objective and highly-constrained settings \citep{chan2018trade,gorissen2013mixed},  including cancer treatment \citep{corletto2003inverse,chan2014generalized, chan2022inverse,ajayi2022objective}, patient care delivery systems \citep{chan2021inverse}, and the diet recommendation problem \citep{ghobadi2018robust,shahmoradi2021quantile, aswani2019behavioral}. 
The diet recommendation problem is further discussed in Section \ref{Section:Application}. For a more comprehensive survey on inverse optimization theory and applications, we refer the readers to the study by \cite{chan2021inversesurvey}.

The inverse learning framework is most closely related to the works of \cite{chan2019inverse} and \cite{shahmoradi2021quantile} in the literature. The work by \cite{chan2019inverse} discusses inverse linear optimization models with noisy data. They employ the recover-then-optimize approach to find the forward optimal solutions which provides little guarantee about compatibility of the final optimal solutions with the given observations, especially when the optimal decision space contains infinitely many points.  \cite{shahmoradi2021quantile}, alternatively, focus on finding extreme-point optimal solutions to address the instability of optimal forward solutions in the recover-then-optimize approach while controlling for outliers. This approach is computationally more expensive and relies on user-given parameters to find the solutions. Additionally, extreme points may bind a high number of constraints which can make the model more goal-driven in the \ObservationConstraintTradeoff. 

\section{Inverse Learning Methodology} \label{Section:Methodology}

The inverse learning framework focuses on recovering an optimal solution and an  underlying cost vector to a partially known linear forward optimization based on a given set of (feasible or infeasible) observations and an optional given set of active constraint priorities (goals). 
To set up the framework, we first detail the forward optimization problem that will be used in both Sections \ref{Section:Methodology} and \ref{sec:AIL}. We then introduce and discuss the inverse learning model. We leave discussions about the \ObservationConstraintTradeoff\ and the corresponding goal-integrated inverse learning to Section \ref{sec:AIL}.  

\subsection{Forward Optimization Model} \label{sec:2.1}
We consider a general forward linear optimization ($\FO$) problem that its feasible region ($\Omega$) is known, non-empty, full dimensional, and free of redundant constraints, but its cost vector ($\bc$) is unknown. This feasible region can be characterized by a set of half-spaces for a linear $\FO$. In this work, and without loss of generality, we divide the set of half-spaces that shape $\Omega$ based on their role and importance in the underlying problem and application. While the exact number of subsets and their importance may be dictated by the problem context, to introduce the methodology, we consider three general constraint types of ({\it i}\,) \emph{\desirable\ constraints} ($\mathcal{R} \neq \emptyset$): a set of constraints that are considered important by the decision-maker, e.g., nutritional constraints in a diet recommendation problem. 
These constraints are analogous to conventional constraints in optimization frameworks and are considered as the default constraint type in this framework too; ({\it ii}\,) \emph{\pref\ constraints} (
$\mathcal{P} \subseteq \mathcal{R}$): a subset of the \desirable\ constraints (if any) that the decision-maker prefers to have active at optimality. 
For instance, goal-specific nutritional constraints (e.g., maximizing the allowed protein intake) can be classified as a \pref\ constraint if the dieter is actively trying to achieve these goals; 
and ({\it iii}\,) \emph{\objective\ constraints} ($\cT$): a set of constraints (if any) that, although necessary for the construction of $\Omega$ and the problem's well-definedness, is not of interest to the decision-maker/user. For instance, non-negativity constraints on food servings variables can be considered \objective\ constraints. 

We note that different sub-classes of \pref\ constraints can also be assumed, each pertaining to a different level of preference, without requiring changes to the methodology. However, for notation simplicity, we consider all \pref\ constraints (if any) to be equally important in the rest of this paper. 
The distinction between \pref\ and \desirable\ constraints will be further explored in Sections \ref{Section:NumericalEx} and \ref{Section:Application} in order to find meaningfully different optimal solutions for the presented examples. In this section, we primarily focus on \desirable\ and \objective\ constraints. If all constraints are of similar importance or there is no knowledge of the significance of the constraints, then we assume $\mathcal{P} = \cT = \emptyset$, and all constraints are contained in the set of relevant constraints $\mathcal{R}$. We formulate the forward optimization model as follows. 
\begin{subequations} \label{FO}
\begin{align}
\FO({\bc,\Omega}): \underset{\bx}{\text{maximize}} & \quad \bc' \bx  \label{eq:FO_obj}\\
\text{subject to} 
& \quad \bA \bx \geq \bb, \label{eq:FO_relevantConst}\\
& \quad \bG \bx \geq \bh,  \label{eq:FO_trivialConst}\\
& \quad \bx \in \mathbb{R}^n.
\end{align}
\end{subequations}
The vector $\bc$ in the cost function \eqref{eq:FO_obj} is assumed to be unknown. The feasible set $\Omega$ is known and characterized as $\Omega =  \left\{\bx \in \mathbb{R}^n| \bA \bx \geq \bb, \bG \bx \geq \bh \right\},$ where  $\bA\in \mathbb{R}^{\left| \mathcal{R}\right|\times n}$, $\bb\in \mathbb{R}^{\left| \mathcal{R}\right|}$ represent the \desirable\ (conventional) constraints $(\left| \mathcal{R}\right| \geq 1)$, and $\bG\in \mathbb{R}^{\left| \cT\right|\times n}$, $\bh\in \mathbb{R}^{\left| \cT\right|}$ the \objective\ constraints ($\cT$). 
Since \pref\ constraints ($\cP$) are a subset of the  \desirable\ constraints $\cR$, for the cases that $\cP \neq \emptyset$, 
we can write 
$\bA = \left [\bA_\cP; \,\bA_{\cR\setminus\cP} \right]$
and $\bb = \left [ \bb_\cP; \,\bb_{\cR\setminus\cP} \right ]$, 
where $\bA_\cP$ and $\bb_\cP$ correspond to the \pref\ constraints and $\bA_{\cR\setminus\cP}$ and $\bb_{\cR\setminus\cP}$ are the not-\pref\ \desirable\ constraints. 
For ease of notation, we define 
$\cJ = \left \{ 1,\dots,\left| \mathcal{R}\right| \right \}$ and 
$\cJ_T = \left \{ 1,\dots,\left| \cT\right| \right \}$ 
to be the set of indices of constraints in $\cR$ and $\cT$, respectively. 
We denote $\Omega^{opt} (\bc)$ as the set of optimal solutions to $\FO(\bc,\Omega)$, and define 
$\Omega^{opt}=\bigcup_{c\neq 0} \Omega^{opt}(\bc)$. For a linear optimization problem, $\Omega^{opt}$ represents the boundary of the closed polyhedral feasible set $\Omega$. 

\subsection{Inverse Learning Model} \label{sec:Methods:mio}

In the inverse setting, our goal is to learn an optimal solution $\bz \in \Omega^{opt}$ and a cost vector $\bc \in \mathbb{R}^n$ of the forward problem \eqref{FO} based on a given set of observation(s), $\bX = \left \{ \bx^k \in \mathbb{R}^{n} : k \in \cK \right \}$, $\cK = \left\{ 1,\hdots, K\right\}$
 which may be noised by error or bounded rationality and not necessarily optimal or feasible for $\FO$. 
Let $\by \in \mathbb{R}_{\geq 0}^{\left| \mathcal{R}\right|}$ and $\bu \in \mathbb{R}_{\geq 0}^{\left| \cT\right|}$ be the corresponding dual variables for the \desirable\ constraints $\cR$ \eqref{eq:FO_relevantConst} and the \objective\ constraints $\cT$ \eqref{eq:FO_trivialConst}, respectively. A solution $\bz \in \Omega$ to $\FO(\bc,\Omega)$ is contained in $\Omega^{opt}$ if and only if a primal-dual feasible pair ($\bz, (\by, \bu)$) exists such that {$\bc' \bz = \by' \bb + \bu' \bh$} and $\bc \neq \bzero$. 
We formulate our \emph{Inverse Learning} model ($\IO$) as follows. 
\begin{subequations} \label{IO}
\begin{align}
\IO({\bX, \Omega}): \underset{\bc, \by, \bu , \bE, \bz}{\text{minimize}} & \quad \cD(\bE)\\
\text{subject to} 
& \quad \bA \bz \geq \bb,   \label{IOPrimalFeasblility1}\\ 
& \quad \bG \bz \geq \bh,   \label{IOPrimalFeasblility2}\\
& \quad \bc' \bz = \bb' \by ,   \label{IOStrongDual}\\ 
& \quad \bz = (\bx^k - \bepsilon^k ),   \quad \forall k \in \mathcal{K}\label{IOOnePoint}\\
& \quad  \bA ' \by  =\bc, \label{SIODualFeas1}\\ 
& \quad {\bu = \bzero}, \label{SIODualFeas2}\\
& \quad \bone_{|\cR|\times 1}' \by = 1,   \label{IORegularization}\\
& \quad \bc, \bz \in \mathbb{R}^n, \quad \by \in \mathbb{R}_{\geq 0}^{\left| \mathcal{R}\right|}, \quad \bu \in \mathbb{R}_{\geq 0}^{\left| \cT\right|}, \quad \bE \in \mathbb{R}^{n\times K}. \label{IOrange}
\end{align}
\end{subequations}
Constraints \eqref{IOPrimalFeasblility1} and \eqref{IOPrimalFeasblility2} enforce primal feasibility for the learned solution $\bz$, as defined by constraint \eqref{IOOnePoint}. 
Constraints \eqref{SIODualFeas1}--\eqref{SIODualFeas2} and \eqref{IOrange} are dual feasibility for relevant and trivial constraints. 
The dual variables corresponding to \desirable\ constraints ($y_{j\in \cJ}$) can assume non-zero values to minimize any influence of the trivial constraints $\cT$ (if any) on the recovered cost vector. 
Constraint \eqref{IOStrongDual} is the strong duality constraint for relevant constraints so that the recovered $\bc$ is not influenced by trivial constraints, therefore the zero term $\bu' \bh$ is dropped based on  constraint \eqref{SIODualFeas2}. 
The matrix of observations' perturbations is denoted by $\bE \in \mathbb{R}^{n \times K}$, where $\epsilon^k$ is the $k^\text{th}$ column of $\bE$  for $k \in \cK$. 
The objective function, $\cD\geq 0$, is a metric of choice, typically a distance function. For example, the sum of norm distances $\cD(\bE) = \sum_{k=1}^{K}{\left \| \epsilon^k \right \|}_l$ or the maximum distance to the observations $\cD(\bE) = max\left \{{\left \| \epsilon^k \right \|}_l : k \in \mathcal{K}\right \}$ to encourage $\bz$ to be uniformly distanced from the observations. 
$\IO$ is feasible for any set of observations, as Proposition \ref{prop:IO_Feas} shows by constructing a feasible solution. 
\begin{proposition}\label{prop:IO_Feas}
For any given $\bX \in \mathbb{R}^{n \times K}$, the solution $(\ba'_j, \besmall _j, \bzero, \left \{ \epsilon_0^1, \cdots, \epsilon_0^K\right \},\bz_0)$ is feasible for $\IO({\bX, \Omega})$, where $\ba_j$ is the $j^{th}$ row of $\bA$, $\epsilon_0^k = \bx^k - \bz_0, \forall k \in \mathcal{K}$ and $j\in\cJ$ corresponds to a constraint in $\cR$ for which $\ba_j \bz_0 = b_j$.
\end{proposition}

One difference between $\IO$ and the existing models lies in the distinction between relevant and trivial constraints, where trivial constraints are not of interest to be active at optimality (constraint \eqref{SIODualFeas2}). Note that this constraint does not mean that \objective\ constraints will not be binding at all for $\bz$, and instead, ensures that the cost vector can be written as a conic combination of only the rows of the relevant matrix, $\bA$. Another notable difference is that a linear constraint, \eqref{IORegularization}, is used to normalize $\bc$ and avoids trivial solutions for the cost vector. While non-convex constraints such as $\left \| \bc \right \|_L=1$ are generally considered in the literature \citep{chan2019inverse,shahmoradi2021quantile,babier2021ensemble}, Proposition  \ref{prop:SIO_regularization} shows that in $\IO$ 
the linear constraint \eqref{IORegularization} is sufficient to rule out trivial solutions such as $\bc = \bzero$. 
This formulation reduces the complexity of the inverse model by replacing the non-convex cost vector regularization constraints with the linear one.
\begin{proposition} \label{prop:SIO_regularization}
For any feasible solution $(\hat{\bc}, \hat{\by}, \hat{\bu}, \hat{\bE}, \hat{\bz})$ for $\IO$, we have $\hat{\bc} \neq \bzero$.
\end{proposition}
Proposition \ref{prop:SIO_regularization} is shown by proving that any solution with $\bc = \bzero$ results in an implicit inequality for $\Omega$ which contradicts $\Omega$ being full dimensional.
%
We note that a desirable characteristic of $\IO$ is that for a feasible solution $(\bc, \by, \bu , \bE, \bz)$ for $\IL$, the learned solution $\bz$ is optimal for $\FO(\bc,\Omega)$, 
and hence, the inverse learning approach extends the existing inverse optimization literature by learning an optimal solution to $\FO$ directly, 
as Theorem \ref{Theorem_IL} illustrates. 

\begin{theorem} \label{Theorem_IL}
If $(\bc^*, \by^*, \bar{\by}^*, \bE^*, \bz^*)$ is optimal for $\IO$, then $\bz^*  \in \Omega^{opt}$.
\end{theorem}

Theorem \ref{Theorem_IL} illustrates that the learned solution through the inverse model is optimal  for the forward model and removes the need to solve it afterwards. 
This distinction is important since solving $\FO(\bc^*,\Omega)$ does not necessarily find the closest optimal solution to $\bX$ if $dim(\Omega^{opt}(\bc^*)) \geq 1$, and is known to be unstable as highlighted in \cite{shahmoradi2021quantile}.  Additionally, since our inverse formulation distinguishes between \desirable\ $\cR$ and \objective\ $\cT$ constraints, the solution $\bz$ can prioritize the relevant constraints to bind at optimality with a recovered $\bc$ that is a conic combination of them. 
As such, the \objective\ constraints that are deemed unimportant by the user have a diminished role. 

Figure \ref{fig:ILvsIO_comparison} illustrates some of the characteristics of the $\IO$ model. 
This schematic figure showcases a given feasible region with two observations (gray dots). 
Employing a 2-norm distance metric, $\cD(\bE) = \sum_{k=1}^{K}{\left \| \epsilon^k \right \|}_2$, we recover the underlying $\FO$ solution and cost vector using classical inverse optimization (black ink) and the $\IO$ model (green ink) for two cases. 
In the first case (Figure \ref{fig:ILvsIO_comparison}(a)), we assume all constraints are conventional, or using our terminology they are all relevant (dashed gray lines), and no expert information or additional knowledge was conveyed about any of the constraints. Solving this problem with both classical inverse optimization and $\IO$ results in the same cost vector. In this example, the main difference between the models is that $\IO$ finds an optimal solution ($\bx_1^*$) that is as close as possible to the observations. The classical inverse optimization returns an extreme point as the optimal solution ($\bar{\bx}$), which is situated further from the observations. 
In the second case (Figure \ref{fig:ILvsIO_comparison}(b)), we assume additional knowledge is available that two of the constraints are indeed trivial (dotted red lines), meaning that they are not important to the user to bind at optimality. For instance, these constraints may be in place to keep the problem well-defined or represent resources that the user has a surplus of and is not concerned with binding them at optimality. The solution found by the classical inverse optimization will not change. The $\IO$, however, can use this additional information to find an optimal solution ($\bx^*_2$) that binds at least one of the relevant constraints. In this example, the corresponding cost vector ($\bc^*_2$) also changes to remain in the cone of relevant constraints instead of being defined solely by trivial constraints. 
\iftrue
\begin{figure}
\begin{center}
\begin{tikzpicture}[scale=0.55] 
  bluenode/.style={shape=rectangle, draw=cyan, line width=2, fill= cyan},
  greennode/.style={shape=bar, draw=gray, line width=2, fill= gray},
  graynode/.style={shape=circle, draw=gray, line width=2, fill= gray},
  blacknode/.style={shape=circle, draw=black, line width=2, fill= black},
  40greennode/.style={shape=circle, draw=black!40!green, line width=2, fill= black!40!green},
  rednode/.style={shape=rectangle, draw=red, line width=2, fill= red}
]

\fill  [black!3!white] (5,1) -- (7,3) -- (7,5)-- (5,7) -- (3,7) -- (1,5) --cycle;
\draw [line width=1pt][gray][costumdashed] (5,1) -- (7,3); 
\draw [line width=1pt][gray][costumdashed] (7,3) -- (7,5);
\draw [line width=1pt][gray][costumdashed] (7,5) -- (5,7);
\draw [line width=1pt][gray][costumdashed] (5,7) -- (2.95,7);
\draw [line width=1pt][gray][costumdashed] (3,7) -- (1,5);
\draw [line width=1pt][gray][costumdashed] (1,5) -- (5,1);
\draw [thick][dashed][black!40!green] (4.75, 4.75) -- (6,6);
\draw [thick][dashed][black!40!green] (4.5, 4.5) -- (6,6);
\draw [thick][dashed][black] (4.75, 4.75) -- (5,7);
\draw [thick][dashed][black] (4.5, 4.5) -- (5,7);
\draw[fill][gray] (4.75, 4.75) circle [radius=0.09];
\draw[fill][gray] (4.5, 4.5) circle [radius=0.09];
\draw[fill][black!40!green] (6, 6) circle [radius=0.15];
\node[label={below,text=black!40!green,label distance=0.1cm}: {$\bx_1^*$}] at (6, 6)  {};
\draw[fill] (5, 7) circle [radius=0.15];
\node[label= {$\Bar{\bx}$}] at (5, 7) {};
\draw[fill][white] (5, 1) circle [radius=0.1];
\node[label={below,text=white}: {$\bx_\CIO$}] at (5, 1) {};

\draw [thick] (5.5,7.5) -- (6.5,6.5);
\draw[->] (6,7)  -- (6.5,7.5) node[right]{$\Bar{\bc}$};
\draw [thick][black!40!green] (6.5,6.5) -- (7.5,5.5);
\draw[->][black!40!green] (7,6)  -- (7.5,6.5) node[right]{$\bc_1^*$};

\node[label={}: {(a) }] at (4,-0.5)  {};
\end{tikzpicture}
\qquad \qquad  
\begin{tikzpicture}[scale=0.55][
  bluenode/.style={shape=rectangle, draw=cyan, line width=2, fill= cyan},
  greennode/.style={shape=bar, draw=gray, line width=2, fill= gray},
  graynode/.style={shape=circle, draw=gray, line width=2, fill= gray},
  blacknode/.style={shape=circle, draw=black, line width=2, fill= black},
  40greennode/.style={shape=circle, draw=black!40!green, line width=2, fill= black!40!green},
  rednode/.style={shape=rectangle, draw=red, line width=2, fill= red}
]

\fill  [black!3!white] (5,1) -- (7,3) -- (7,5)-- (5,7) -- (3,7) -- (1,5) --cycle;
\draw [line width=1pt][gray][costumdashed] (5,1) -- (7,3);
\draw [line width=1pt][gray][costumdashed] (7,3) -- (7,5);
\draw [line width=1pt][red][very loosely dotted] (7,5) -- (5,7);
\draw [line width=1pt][red][very loosely dotted] (5,7) -- (2.95,7);
\draw [line width=1pt][gray][costumdashed] (3,7) -- (1,5);
\draw [line width=1pt][gray][costumdashed] (1,5) -- (5,1);

\draw [thick][dashed][black!40!green] (4.75, 4.75) -- (7, 5);
\draw [thick][dashed][black!40!green] (4.5, 4.5) -- (7, 5);
\draw [thick][dashed][black] (4.75, 4.75) -- (5,7);
\draw [thick][dashed][black] (4.5, 4.5) -- (5,7);
\draw[fill][gray] (4.75, 4.75) circle [radius=0.09];
\draw[fill][gray] (4.5, 4.5) circle [radius=0.09];
\draw[fill][black!40!green] (7, 5) circle [radius=0.15];
\node[label={above,text=black!40!green,label distance=0.01cm}: {$\bx_2^*$}] at (7, 5)  {};
\draw[fill] (5, 7) circle [radius=0.15];
\node[label= {$\Bar{\bx}$}] at (5, 7) {};
\draw[fill][white] (5, 1) circle [radius=0.1];
\node[label={below,text=white}: {$\bx_\CIO$}] at (5, 1) {};
\draw [thick] (5.5,7.5) -- (6.5,6.5);
\draw[->] (6,7)  -- (6.5,7.5) node[right]{$\Bar{\bc}$};
\draw [thick][black!40!green] (7.5,5) -- (7.5,4);
\draw[->][black!40!green] (7.5,4.5)  -- (8,4.5) node[right]{$\bc_2^*$};
\node[label={}: {(b)}] at (4,-0.5)  {};
\end{tikzpicture}
\begin{tikzpicture}[scale=0.25][] 
\footnotesize
\draw[fill][gray] (9.5,15) circle [radius=0.2];  \node[label={right,label distance=0.1cm}: {Observation}] at (10.5,15)  {};
\draw[fill][black!30!green] (9.5,13.6) circle [radius=0.3]; \node[label={right,label distance=0.1cm}: {$\IO$ Output}] at (10.5,13.6)  {};
\draw[fill][black] (9.5,12.2) circle [radius=0.3]; \node[label={right,label distance=0.1cm}: {$\CIO$ Output}] at (10.5,12.2)  {};
\draw [line width=1pt][red][very loosely dotted] (8.4,10.8) -- (11.4,10.8); \node[label={right,label distance=0.1cm}: {\Objective\ Constraint}] at (10.5,10.8)  {};
\draw [line width=1pt][gray][costumdashed] (8.4,9.4) -- (11.4,9.4); \node[label={right,label distance=0.1cm}: {\Desirable\ Constraint}] at (10.5,9.4)  {};
\draw [line width=4pt][white] (8.4,4) -- (9,4);
\draw [line width=0.5][white] (8,16) -- (22,16);
\draw [line width=0.5][white] (22,16) -- (22,8);
\draw [line width=0.5][white] (22,8) -- (8,8);
\draw [line width=0.5][white] (8,16) -- (8,8);
\end{tikzpicture}
\caption{The figures show a given feasible region (shaded gray) and two observations (gray dots). Using a 2-norm distance metric, classical inverse optimization recovers $\bar{c}$ and $\bar{x}$ as the cost vector and optimal solution of $\FO$. The inverse learning approach (green ink), however, finds a closer optimal solution ($\bx^*_1$) in Figure (a), even though the cost vector $\bc^*_1$ remains similar to that recovered by inverse optimization. In Figure (b), additional information is provided that two of the constraints are trivial (red dotted lines). In this case, the inverse learning approach adjusts and recovers the closest point that binds at least one relevant constraint (dashed gray). }
\label{fig:ILvsIO_comparison}
\end{center}
\end{figure}
\fi

\medskip
{\bf Solving Inverse Learning Model.} To solve the $\IO$ models, we note that constraint \eqref{IOStrongDual} is bilinear and non-convex, analogous to other inverse models in the literature. 
A similar approach to \cite{chan2019inverse} can be adapted for the $\IO$ framework and relevant constraints to solve $\IO$ by solving a series of $m_{1}$ linearly constrained problems. For any feasible solution $(\bc, \by, \bar{\by} , \bE, \bz)$, we show that there exists $j \in \cJ$ such that $\ba_j \bz = b_j$, where $\ba_j$ is the $j^{th}$ row of the matrix $\bA$. Therefore, $\bz \in \Omega^{opt}(\ba_j), ~\forall j \in J$. Considering solutions that bind at least one \desirable\ constraint allows for removing duality and complementary slackness constraints from \IO. As such, there exists an optimal solution $(\bc, \by, \bar{\by} , \bE, \bz)$ for $\IO$ where $\bc$ is orthogonal to at least one of the \desirable\ constraints of $\FO$, as formalized in Proposition \ref{theorem1}. 
\begin{proposition}\label{theorem1}
Let $\bc = {\ba'_j}$ and $\by = {e_j}$ where $j \in \cJ$ designates the $j^{th}$ \desirable\ constraint and $e_j$ is the $j^{th}$ unit vector. Formulation \eqref{MFD} finds $\bz^j \in \Omega$ which binds the j$^\text{th}$ constraint in $\cR$ and minimizes $\cD$. Let $j_{min} = argmin_{j \in \cJ} \left\{ D(E_j) \right\}$  and $(\bE^*$, $\bz^*)$ be the optimal solution of $\IO_{j_{min}}$, then if the optimal objective value of $\IO$ is $\cD^*$, we have $\cD_{min} = \cD^*$ and  $(\ba_{j_{min}}$, $e_{j_{min}}$, $\bzero$, $\bE^*$, $\bz^*)$ is optimal for $\IO$.
\begin{subequations} \label{MFD}
\begin{align}
\IO_j({\bX,\Omega}): 
\underset{\bz, \bE_j}{\text{minimize}} & \quad \cD(\bE_j)
\\
\text{subject to} 
& \quad \ba_j \bz^j = b_j,  \label{ILjStrongDual}\\ 
& \quad \bA \bz^j \geq \bb,   \label{ILjPrimal Feasblility1}\\ 
& \quad \bG \bz^j \geq \bh,   \label{ILjPrimal Feasblility2}\\
& \quad \bz^j = (\bx^k - \bepsilon^k_j ),  \quad \forall k \in \mathcal{K}\label{ILjOnePoint} \\
& \quad \bz^j \in \mathbb{R}^n,  \quad \bE_j \in \mathbb{R}^{n\times K}. \label{ILjrange}
\end{align}
\end{subequations}
\end{proposition}

Proposition \ref{theorem1} demonstrates that since all learned solutions from $\IO$ bind at least one \desirable\ constraint, it suffices to solve simplified linearly constrained optimization models that find $\FO$ solutions that bind specified \desirable\ constraint while minimizing the metric $\cD$. As such, an optimal solution for $\FO$ (which is a part of an optimal solution of $\IL$) can be achieved by solving such optimization problems for all \desirable\ constraints and choosing the solution corresponding to the minimum $\cD$ among the all obtained solutions. Furthermore, formulation \eqref{MFD} is a linearly constrained problem which aids its computational complexity. 

\section{Goal-integrated Inverse Learning} \label{sec:AIL}

Section \ref{Section:Methodology} focuses on learning forward optimal solutions that are as close as possible to the given observations while distinguishing between relevant and trivial constraints. 
For a polyhedral feasible set $\Omega$, binding fewer constraints translates to the ability to be closer to the observations. Hence, relying more on the observations in the learning process. In contrast, binding more constraints represents trusting the known information on the feasible set and relying less on the observations. 
Therefore, there is an inherent tradeoff between learning from the observations and the feasible set. As an example, given known dietary decisions in a diet recommendation problem, $\IO$ recommends diets that optimize a cost function while minimizing a given distance metric with regard to the observations. Conversely, if we aim to recommend highly nutritious diets, then the learned solution may need to bind as many nutritional constraints as possible, and as a result, shift away from the observations, and hence, become less personalized to the patient's past behavior. 
We refer to this interplay between the goal- or observation-driven solutions as the \emph{\ObservationConstraintTradeoff} and note that binding more constraints may come at the cost of increasing the distance from the observations, i.e., increasing the value of $\IO$'s objective metric $\cD$ at optimality. 

In this section, we explore embedding constraint goals into the inverse models and controlling the number and the type of binding constraints (Section \ref{sec:MGIL}). We analyze the  \ObservationConstraintTradeoff\ and utilize it to provide a series of adaptive models that are capable of maintaining past goals as new ones are added, hence, lending themselves to applied settings where gradual improvements are favored (Section \ref{Sec:GIL_variants}). 


\subsection{The Goal-integrated Inverse Learning Model} \label{sec:MGIL}

We develop the \emph{Goal-integrated Inverse Learning} ($\MGIL$) to balance the \ObservationConstraintTradeoff. The model provides a mechanism to control the number of active relevant constraints while prioritizing binding those constraints that are preferred ($\mathcal{P}\subset \mathcal{R}$).   
To regulate the number of binding \desirable\ constraints, we introduce a set of binary variables, $v_j$ for all $j \in \cJ$, and an input parameter $1 \leq r \leq n$. 
%
Let $\bP \subseteq \cJ$ be the set of indices of \pref\ \desirable\ constraints, $\bv$ the vector of binary variables $v_j$, $j\in \cJ$, and $\bE$ the perturbation matrix as defined in Section \ref{sec:Methods:mio}. The $\MGIL$ model can be written as follows.
%
\begin{subequations} \label{MGIL}
\begin{align}
\MGIL({\bX, \Omega, r,\bP}): 
\underset{\mathbf{v}, \bE, \bz}{\text{minimize}} & \quad  \omega \cD(\bE) - (1-\omega)\sum_{s \in \bP} v_s 
\\
\text{subject to} 
& \quad \bb \leq \bA  \bz \leq \bb + M (1-\bv), \label{GILStrongDual}\\ 
& \quad \bG \bz \geq \bh, \label{GILPrimal Feasblility2}\\
& \quad \bz = (\bx^k - \bepsilon^k ),   \quad \forall k \in \mathcal{K}\label{GILOnePoint}\\
& \quad \sum_{j \in \cJ} v_j =  r,    \label{GILCornerPoint}\\ 
& \quad  \bz \in \mathbb{R}^n,  \quad \bE \in \mathbb{R}^{n\times K}, \quad \bv \in \{0,1\}^{|\cJ|}. \label{GILrange}
\end{align}
\end{subequations}
Constraints \eqref{GILStrongDual} and \eqref{GILCornerPoint} ensure that exactly $r$ \desirable\ constraints are binding for $\bz$. If no knowledge of the exact value for $r$ is available, we can either assume $r=1$ or convert constraint \eqref{GILCornerPoint} to a Lagrangian multiplier to maximize the number of binding \desirable\ constraints. 
The objective function of $\MGIL$ forces the learned solution $\bz$ to minimize $\cD$ while binding as many \pref\ constraints as possible. The parameter $\omega$ is a user-defined weight that adjusts the \ObservationConstraintTradeoff\ to the desirable balance between the distance to the observations and binding additional \pref\ constraints. 
Figure \ref{fig:IL_general_goals} illustrates the \ObservationConstraintTradeoff\ 
in a schematic two-dimensional example with a single observation (gray dot). The trivial constraints are shown in red dots, relevant constraints in dashed gray, and the \pref\ constraint in the solid black line in the figures, and a 2-norm metric is used as the distance function. The preferred constraint is not highlighted for binding (i.e., $\omega=1$) in Figure \ref{fig:IL_general_goals}(a). Instead, the model is solved to bind one and then two \desirable\ constraints, $\bz|_{r=1}$ and $\bz|_{r=2}$, respectively. As the figure shows, as the value of $r$ increases, the distance between observations and the learned solution may also increase as the solutions rely more on constraints information as opposed to the observations. Figure \ref{fig:IL_general_goals}(b) depicts the same model but when the preferred constraint is engaged. The model aims to bind the preferred constraint when possible. 
This preference may increase the distance from observations depending on the location of the preferred constraints. 
\begin{figure}[t]
\begin{center}
\begin{tikzpicture}[scale=0.5][
  bluenode/.style={shape=rectangle, draw=cyan, line width=2, fill= cyan},
  greennode/.style={shape=bar, draw=gray, line width=2, fill= gray},
  graynode/.style={shape=circle, draw=gray, line width=2, fill= gray},
  blacknode/.style={shape=circle, draw=black, line width=2, fill= black},
  40greennode/.style={shape=circle, draw=black!40!green, line width=2, fill= black!40!green},
  rednode/.style={shape=rectangle, draw=red, line width=2, fill= red}
]

\fill  [black!5!white] (7,1) -- (7,4) -- (4,7)-- (5,7) -- (1,7) -- (1,1) -- cycle;
\draw [line width=1pt][gray][costumdashed] (7,1) -- (7,4);
\draw [line width=1pt][gray][costumdashed] (7,4) -- (4,7);
\draw [line width=1pt][black] (4,7) -- (1,7);
\draw [line width=1pt][red][very loosely dotted] (1,7) -- (1,1);
\draw [line width=1pt][red][very loosely dotted](1,1) -- (7,1);

\draw [very thick][dashed][gray] (4.5,4.5) -- (5.5,5.5);
\draw[fill][black!30!green] (5.5,5.5) circle [radius=0.15];
\node[label=right: {$\bz|_{r=1}$}] at (5.5,5.5) {};

\draw [very thick][dashed][gray] (4.5,4.5) -- (7,4);
\draw[fill][black!30!green] (7,4) circle [radius=0.15];
\node[label=right: {$\bz|_{r=2}$}] at (7,4) {};

\draw[fill][gray] (4.5,4.5) circle [radius=0.1];

\node[label={}: {\footnotesize  (a) Constraint-Observation tradeoff }] at (4,-0.5)  {};
\end{tikzpicture}
\qquad
\begin{tikzpicture}[scale=0.5][
  bluenode/.style={shape=rectangle, draw=cyan, line width=2, fill= cyan},
  greennode/.style={shape=bar, draw=gray, line width=2, fill= gray},
  graynode/.style={shape=circle, draw=gray, line width=2, fill= gray},
  blacknode/.style={shape=circle, draw=black, line width=2, fill= black},
  40greennode/.style={shape=circle, draw=black!40!green, line width=2, fill= black!40!green},
  rednode/.style={shape=rectangle, draw=red, line width=2, fill= red}
]
\fill  [black!5!white] (7,1) -- (7,4) -- (4,7)-- (5,7) -- (1,7) -- (1,1) -- cycle;
\draw [line width=1pt][gray][costumdashed] (7,1) -- (7,4);
\draw [line width=1pt][gray][costumdashed] (7,4) -- (4,7);
\draw [line width=1pt][black] (4,7) -- (1,7);
\draw [line width=1pt][red][very loosely dotted] (1,7) -- (1,1);
\draw [line width=1pt][red][very loosely dotted](1,1) -- (7,1);

\draw [very thick][dashed][gray] (4.5,4.5) -- (4,7);
\draw[fill][black!30!green] (4,7) circle [radius=0.15];
\node[label=right: {$\bz$}] at (4,7) {};

\draw[fill][gray] (4.5,4.5) circle [radius=0.1];

\node[label={}: {\footnotesize (b) Binding Preferred Constraints}] at (4,-0.5)  {};
\end{tikzpicture}
\begin{tikzpicture}[scale=0.25][] 
\footnotesize
\draw[fill][gray] (9.5,15) circle [radius=0.2];  \node[label={right,label distance=0.1cm}: {Observation}] at (10.5,15)  {};
\draw[fill][black!30!green] (9.5,13.6) circle [radius=0.3]; \node[label={right,label distance=0.1cm}: {$\MGIL$ Output}] at (10.5,13.6)  {};
\draw [line width=1pt][red][very loosely dotted] (8.4,12.2) -- (11.4,12.2); \node[label={right,label distance=0.1cm}: {\Objective\ Constraint}] at (10.5,12.2)  {};
\draw [line width=1pt][gray][costumdashed] (8.4,10.8) -- (11.4,10.8); \node[label={right,label distance=0.1cm}: {\Desirable\ Constraint}] at (10.5,10.8)  {};
\draw [line width=1pt][black] (8.4,9.4) -- (11.1,9.4); \node[label={right,label distance=0.1cm}: {\Pref\ Constraint}] at (10.5,9.4)  {};
\draw [line width=4pt][white] (8.4,4) -- (9,4);
\draw [line width=0.5][white] (8,16) -- (22,16);
\draw [line width=0.5][white] (22,16) -- (22,8);
\draw [line width=0.5][white] (22,8) -- (8,8);
\draw [line width=0.5][white] (8,16) -- (8,8);
\end{tikzpicture}
\caption{
A schematic figure to illustrate the \ObservationConstraintTradeoff\ in $\MGIL$. Figure (a) shows the results for $r=1,2$ when the \pref\ constraint is not activated ($\omega = 1$). 
The learned solution $\bz|_{r=1}$ binds one \desirable\ constraint (dashed gray lines) and $\bz|_{r=2}$ binds two, although the distance to observations increases. 
Figure (b) illustrates the model with an activated \pref\ constraint ($\omega \neq 1$). The learned solution prioritizes binding the \pref\ constraint (solid black lines) when possible.
}
\label{fig:IL_general_goals}
\end{center}
\end{figure}

Given that the feasibility of the $\MGIL$ model may depend on the input parameter $r$, Proposition \ref{MGIL_feas} discusses the feasibility of the model under various conditions.
\begin{proposition} \label{MGIL_feas} 
We have the following for $\MGIL$: 
\begin{enumerate}[(a)]
    \item  $\MGIL(\bX,\Omega,r,\bP)$ is feasible for $r = 1$.
    \item $\MGIL(\bX,\Omega,r,\bP)$ is feasible for any $\forall r \in \left \{ 1,...,n \right \}$ when $\cT = \emptyset$.
    \item  $\MGIL(\bX,\Omega,r,\bP)$ is feasible for $r \in \left \{ 1,...,r_0 \right \}$ if $\Omega$ has at least one face $\cF$ with $r_0$ binding \desirable\ constraints. In particular, $\MGIL$ is feasible for all $r \in \left \{ 1,...,n \right \}$ if $\Omega$ has at least one extreme point with $n$ binding \desirable\ constraints.
    \item If $\cT = \emptyset$ and $1 \leq r_1 \leq r_2 \leq n$, then $\cD^*_{r_1} \leq \cD^*_{r_2}$, where $\cD^*_{r_1}$ and $\cD^*_{r_2}$ are the optimal objective terms $\cD(\bE)$ in $\MGIL(\bX,\Omega,r_1,\bP)$ and $\MGIL(\bX,\Omega,r_2,\bP)$, respectively. 
\end{enumerate}
\end{proposition}
As Proposition \ref{MGIL_feas} shows, the $\MGIL$ model is always feasible for $r=1$. In general, the feasibility of $\MGIL$ is contingent on the existence of a face of $\Omega$ with $r$ \desirable\ binding constraints. Although showing that such faces exist for different values of $r$ is non-trivial \citep{civril2009selecting}, Proposition \ref{MGIL_feas}(b) asserts that in the absence of \objective\ constraints where all constraints are assumed to be \desirable\ and conventional, $\MGIL$ is always feasible. 
Proposition \ref{MGIL_feas}(d) summarizes the tradeoff between distance to the observations and binding more constraints at optimality for distance-based metrics. 
The \ObservationConstraintTradeoff\ is controlled by the input parameter $r$, and as the value of $r$ increases, so does the value of $\cD$ at optimality. 

We note that $\MGIL$ does not explicitly recover a cost vector. Instead, since one can recover the binding constraints of any $\bz \in \Omega^{opt}$, we can readily find a cost vector $\bc$ such that $\bz$ is optimal for $\FO(\bc,\Omega)$ using previous results \citep{tavasliouglu2018structure}. Theorem \ref{Theorem2} characterizes the set of all such cost vectors. We use the notation 
$\cJ_z \subseteq \cJ$ to characterize the set of all indices of \desirable\ constraints that are binding for $\bz$. 

\begin{theorem} \label{Theorem2}
Let $(\mathbf{v}^*, \bE^*, \bz^*)$ be optimal for $\MGIL (\bX,\Omega,r_0,\bP)$. Then, $\bz^*$ is optimal for $\FO (\bc,\Omega)$ for any $\bc \in cone (\ba_t : t \in \cJ_{z^*})$.
\end{theorem}

Theorem \ref{Theorem2} illustrates two important properties of the goal-integrated inverse  learning model. First, $\MGIL$ indeed learns an optimal solution to $\FO(\bc,\Omega)$, where $\bc$ is contained in the recovered cone. 
This result asserts that $\MGIL$ is capable of solving the inverse problem for recovering the missing parameters and learning an optimal solution. Second, Theorem \ref{Theorem2} shows that learning the optimal solution for $\FO$ is sufficient to characterize a set of inverse optimal cost vectors. The user can then choose their desired cost vector to recover $\FO$, using their method of choice. If no such method exists, we propose the following method to select a single $\bc$ from the recovered cone.  
Let $\bz^* \in \Omega^{opt}$ be the learned solution from $\MGIL(\bX,\Omega,r,\bP)$, then the set of non-trivial $\bc$ which makes $\bz^*$ optimal for $\FO (\bc,\Omega)$ is characterized by Theorem \ref{Theorem2} as
\begin{equation} \label{InverseFeasibleRegion}
 \quad \cQ (\bz^*) = 
\left \{ 
\left.\begin{matrix} \bc \in \mathbb{R}^{n} \end{matrix}\right| 
\bc \in cone (\ba_t : t \in \cJ_{z^*}), \left \| \bc \right \|_L=1
\right \},
\end{equation} 
where a normalization $\|\cdot\|_L$ is added so that each direction of $\bc$ is represented exactly once in $\cQ$. To characterize a specific $\bc \in \cQ$, for instance, the objective value of the observations ($\sum_{k=1}^{K} \bc'\bx^k$) can be maximized by the following model.
\begin{subequations} \label{InverseOptimalc}
\begin{align}
\underset{\bc}{\text{maximize}} & 
\quad \sum_{k=1}^{K} \bc'\bx^k\\
\text{subject to} 
& \quad \bc \in \cQ (\bz^*).
\end{align}
\end{subequations}
We keep the general definition of $\cQ$ for $\MGIL$ as a generalization of $\IO$. $\IO$ is a reduced version of $\MGIL$ when $\bP=\emptyset$, i.e., no preferred constraints are indicated either due to a lack of knowledge or the choice of the user to assume no \pref\ constraint. Theorem \ref{rem:GILlinktoIL} summarizes this result. 
\begin{theorem}\label{rem:GILlinktoIL}
A solution $(\mathbf{v}_{\MGIL}$ $\bE_{\MGIL}$, $\bz_{\MGIL})$ is optimal for $\MGIL(\bX,\Omega,r = 1, \bP = \emptyset)$ if and only if there exists $(\bc_{\IO}$, $\by_{\IO}$, $\bu_{\IO}$, $\bE_{\IO}$, $\bz_{\IO})$ optimal for $\IO(\bX,\Omega)$ such that $\bz_{\IO} = \bz_{\MGIL}$.
\end{theorem}
Theorem \ref{rem:GILlinktoIL} shows that $\MGIL$ generalizes $\IO$ for learning $\bz \in \Omega^{opt}$. Since $\IO$ learns a solution that minimizes the metric $\cD$, Theorem \ref{rem:GILlinktoIL} shows that $\MGIL(\bX,\Omega,r = 1, \bP = \emptyset)$ also finds a solution $\bz \in \Omega^{opt}$ that minimizes $\cD$ when $r=1$ for the given $\bX$, hence the two models become equivalent in solving $\bz$ as Corollary \ref{Corollary 1} shows. 
\begin{corollary} \label{Corollary 1}
Let $(\mathbf{v}^*, \bE^*, \bz^*)$ be optimal for $\MGIL(\bX,\Omega,r = 1, \bP = \emptyset)$, then $\forall \hat{\bc} \in cone (\ba_t : t \in \cJ_{z^*})$ where $\hat{\bc} \neq \bzero$, we have $\hat{\bc} = \bA' \hat{\by}$ where $\hat{\by} \in \mathbb{R}_{\geq 0}^{\left| \mathcal{R}\right|}$ such that $\hat{y}_j = 0, \forall j \notin \cJ_{z^*}$ and  $(\hat{\bc},\hat{\by}, \bzero, \bE^*, \bz^*)$ is optimal for $\IO(\bX,\Omega)$. 
\end{corollary}
Using the expression of Corollary \ref{Corollary 1}, we observe that the learned solution $\bz$ applying $\MGIL$ to a given set of observations $\bX$ does in fact satisfy strong duality conditions, and consequently, the complementary slackness constraints for $\FO(\bc,\Omega)$ for any non-zero $\bc \in cone (\ba_t : t \in \cJ_{z^*})$. 

\medskip
{\bf The Observation-goal tradeoff.} For the remainder of this section, we turn our attention to the interplay between the behavior of $\MGIL$ and different choices for $\bP$. 
Remark \ref{prop:MGILvsGIL} shows that $\MGIL$ learns a solution $\bz \in \Omega^{opt}$ that returns a larger value of the metric $\cD$ when there are \pref\ constraints. In other words, when $\bP\neq\emptyset$, then the solution of $\MGIL$ is no better than when $\bP=\emptyset$ because inclusion of \pref\ constraints results in learning solutions with increased values of the distance metric $\cD$ at optimality. 

Figure \ref{fig:MGIL} shows an example of how learned solutions of $\MGIL$ compare against each other for the case of a 2-norm distance metric. 
Three learned solutions are demonstrated in Figure \ref{fig:MGIL}, where $(\bv_1, \bE_1,\bz_1)$ and $(\bv_2, \bE_2,\bz_2)$ are optimal solutions to $\MGIL(\bX,\Omega,r = 1, \bP = \emptyset)$ and $\MGIL(\bX,\Omega,r = 2, \bP = \emptyset)$, respectively, and $(\bv_3, \bE_3,\bz_3)$ is an optimal solution to $\MGIL(\bX,\Omega,r, \bP \neq \emptyset)$. 
A feasible region with trivial (red dots) and relevant (dashed gray) constraints is considered. For cases where $\bP \neq \emptyset$, the \pref\ constraint is shown with a solid line. The figure demonstrates the \ObservationConstraintTradeoff\ with distance to observations increasing as the number of relevant binding constraints or \pref\ constraints increases, as also highlighted in Remark \ref{prop:MGILvsGIL}. 
\begin{remark}\label{prop:MGILvsGIL}
Let $(\mathbf{v}^1_{\MGIL}, \bE^1_{\MGIL}, \bz^1_{\MGIL})$ and $(\mathbf{v}^2_{\MGIL}, \bE^2_{\MGIL}, \bz^2_{\MGIL})$ be optimal for $\MGIL_1(\bX,\Omega,r,\bP_1 = \emptyset)$ and $\MGIL_2(\bX,\Omega,r, \bP_2 \neq \emptyset)$ respectively. Then, $\cD (\bE^1_{\MGIL}) \leq \cD (\bE^2_{\MGIL})$.
\end{remark}

\iftrue
\begin{figure}
\footnotesize
\begin{center}

\begin{tikzpicture}[scale=0.7][
  bluenode/.style={shape=rectangle, draw=cyan, line width=2, fill= cyan},
  greennode/.style={shape=rectangle, draw=gray, line width=2, fill= gray},
  graynode/.style={shape=circle, draw=gray, line width=2, fill= gray},
  blacknode/.style={shape=circle, draw=black, line width=2, fill= black},
  rednode/.style={shape=rectangle, draw=red, line width=2, fill= red}
]
\fill  [black!3!white] (5,1) -- (7,3) -- (7,5)-- (5,7) -- (3,7) -- (1,5) --cycle;
\draw [line width=1pt][red][very loosely dotted] (5,1) -- (7,3);
\draw [line width=1pt][gray][costumdashed] (7,3) -- (7,5);
\draw [line width=1pt][gray][costumdashed] (7,5) -- (5,7);
\draw [line width=1pt][black](5,7) -- (2.95,7);
\draw [line width=1pt][gray][costumdashed] (3,7) -- (1,5);
\draw [line width=1pt][red][very loosely dotted] (1,5) -- (5,1);

\draw[fill][gray] (5.5, 3.25) circle [radius=0.0667];
\draw[fill][gray] (5.5, 3.75) circle [radius=0.0667];

\draw [thick][dashed][gray] (5.5, 3.25) -- (7, 3.5);
\draw [thick][dashed][gray] (5.5, 3.75) -- (7, 3.5);
\node[label=below: {$\cD(\bE_1)$}] at (6.25, 3.5) {};
\draw [thick][dashed][gray] (5.5, 3.25) -- (7, 5);
\draw [thick][dashed][gray] (5.5, 3.75) -- (7, 5);
\node[label=above: {$\cD(\bE_2)$}] at (6, 4.25) {};
\draw [thick][dashed][gray] (5.5, 3.25) -- (5, 7);
\draw [thick][dashed][gray] (5.5, 3.75) -- (5, 7);
\node[label=left: {$\cD(\bE_3)$}] at (5.3, 5.5) {};

\draw[fill][black!30!green] (7, 3.5) circle [radius=0.1];
\node[label=right: {$\bz|_{r = 1}^{\bP = \emptyset}$}] at (7, 3.5) {};
\draw[fill][black!30!green] (7, 5) circle [radius=0.1];
\node[label=right: {$\bz|_{r = 2}^{\bP = \emptyset}$}] at (7, 5) {};
\draw[fill][black!30!green] (5, 7) circle [radius=0.1];
\node[label= {$\bz|_{r = 1}^{\bP \neq \emptyset}$}] at (5, 7) {};
\end{tikzpicture}
\begin{tikzpicture}[scale=0.3][] 
\footnotesize
\draw[fill][gray] (11.5,15) circle [radius=0.2];  \node[label={right,label distance=0.1cm}: {Observation}] at (12.5,15)  {};
\draw[fill][black!30!green] (11.5,13.6) circle [radius=0.3]; \node[label={right,label distance=0.1cm}: {$\MGIL$ Output}] at (12.5,13.6)  {};
\draw [line width=1pt][red][very loosely dotted] (10.4,12.2) -- (13.4,12.2); \node[label={right,label distance=0.1cm}: {\Objective\ Constraint}] at (12.5,12.2)  {};
\draw [line width=1pt][gray][costumdashed] (10.4,10.8) -- (13.4,10.8); \node[label={right,label distance=0.1cm}: {\Desirable\ Constraint}] at (12.5,10.8)  {};
\draw [line width=1pt][black] (10.4,9.4) -- (13.1,9.4); \node[label={right,label distance=0.1cm}: {\Pref\ Constraint}] at (12.5,9.4)  {};
\draw [line width=4pt][white] (8.4,4) -- (9,4);
\draw [line width=0.5][white] (8,16) -- (22,16);
\draw [line width=0.5][white] (22,16) -- (22,8);
\draw [line width=0.5][white] (22,8) -- (8,8);
\draw [line width=0.5][white] (8,16) -- (8,8);
\end{tikzpicture}
\caption{
The learned optimal solutions $\bz_1$ and $\bz_2$ bind one and two \desirable\ constraints, respectively, and $\bz_3$ binds two \desirable\ constraints with one of them being the sole \pref\ constraint. Although, the distance to observation increases with $\cD(\bE_3) > \cD(\bE_2) > \cD(\bE_1)$. 
} 
\label{fig:MGIL}
\end{center}
\end{figure}
\fi

We illustrated that $\MGIL$ generalizes the inverse optimization problem, however, in the presence of a high number of \desirable\ constraints, different faces of $\Omega$ may bind for different values of $r$. While this effect is in accordance with the \ObservationConstraintTradeoff, it makes comparisons and interpretations of the learned solutions in terms of their binding constraints and slack values more challenging. 
In the next section, we further study $\MGIL$ derivatives that enhance user's ability to compare the learned solutions for different values of $r$. The adaptive model provides a platform to successively increase $r$ to compound more goals while navigating the tradeoff of each additional goal.

\subsection{Adaptive Sequencing of Learned Optimal Solutions} \label{Sec:GIL_variants}

The $\MGIL$ model learns optimal solutions of the forward optimization problem when the user provides one set of relevant ($\cR$), trivial ($\cT$), and \pref\ constraints ($\mathcal{P} \subseteq \cR$). 
In settings where additional solutions are desired or those in which the optimal solution includes a high number of relevant constraints that draws it away from the observations, it is paramount to have a model that can balance the \ObservationConstraintTradeoff\ in a structured manner. In particular, a model that can adapt prior solutions to bind more relevant constraints without losing the already binding constraints can better contrast the tradeoff and accommodate a more informed comparison between these solutions, by keeping them contained in a single face of $\Omega$. The increased control also enables obtaining a sequence of solutions that spans the \ObservationConstraintTradeoff\ space. 
%

The application of the diet recommendation problem is a prime example of a case where comparability of solutions for different values of $r$ is essential. 
For instance, if a learned solution from $\MGIL(\bX,\Omega,1,\bP)$ only binds the maximum protein limit constraint and the learned solution from $\MGIL(\bX,\Omega,2,\bP)$ binds the nutritional constraints for minimum sodium and maximum iron levels, the process of deciding between the two solutions is non-intuitive unless the user has quantified preferences over these nutritional constraints. However, if the learned solution from $\MGIL(\bX,\Omega,2,\bP)$ also binds the maximum protein constraint, the user can directly assess the cost of binding additional constraints in terms of the metric $\cD$. In addition to enabling comparison between different solutions, adapting the $\MGIL$ to keep the previous binding constraints provides a base for the users to traverse the space between observation-driven to goal-driven solutions. The obtained sequence of solutions as $r$ increases provides a path for the users to start as close as possible to the observations and their habits and gradually bind more constraints to reach a more conventional (goal-driven) optimal solution. 

In more detail, let $\bz_0$ be a known solution to $\FO(\bc,\Omega)$ for $\bc \in cone (\ba_t : t \in \cJ_{z_0})$ such that $\cJ_{z_0} \neq \emptyset$. 
Such a solution $\bz_0$ might be available to the user from prior knowledge or might be the result of solving $\MGIL(\bX,\Omega,r_0,\bP)$ for a known set of observations.  
The adaptive $\MGIL$ model, $\MGIL_r$, takes the previous solution, $\bz_r$, as an input and finds a new solution such that at least one additional \desirable\ constraint is appended to the set of previously binding constraints. 
For given $\bz_r$, a new solution $\bz_{r+1}$ is found by shifting $\bz_{r}$ in a direction that binds at least one more \desirable\ constraint with minimal perturbation from $\bX$, if such a shift is possible, while maintaining all previously binding \desirable\ constraints. 
The adaptive $\MGIL_r$ model can be formulated as
\begin{subequations} \label{GGIL}
\begin{align}
\GGIL({\bX, \Omega, \bz_{r}, \bP}): \underset{\mathbf{v}, \bE, \bz_{r+1}}{\text{minimize}} & \quad  \omega \cD(\bE) - (1-\omega)\sum_{s \in \bP} v_s 
\\
\text{subject to} 
& \quad \bb \leq \bA  \bz_{r+1} \leq \bb + M (1-\bv), \label{GGILStrongDual}\\ 
& \quad \bG \bz \geq \bh,   \label{GGILPrimal Feasblility2}\\
& \quad \bz_{r+1} = (\bx^k - \bepsilon^k ),   \quad \forall k \in \mathcal{K}\label{GGILOnePoint}\\
& \quad \sum_{j \in \cJ} v_j \geq  \left | \cJ_{z_i} \right | + 1,    \label{GGILCornerPoint}\\ 
& \quad \cJ_{z_r} = \left \{ j \in \cJ | \ba_j \bz_r = b_j \right \}, \label{GGILBinding} \\
& \quad v_j =  1,  \quad \forall j \in \cJ_{z_r}, \quad \bv \in \{0,1\}^{|\cJ|}.
\label{GGILDependence} 
\end{align}
\end{subequations}

Formulation \eqref{GGIL} inputs a known solution $\bz_r$ that is optimal for $\FO(\bc,\Omega)$ and learns a solution $\bz_{r+1}$ that binds more \desirable\ constraints. 
Given an initial solution $\bz_r$, 
a sequence of solutions $\bz_r,\bz_{r+1}, \hdots, \bz_{L}$ can be achieved where for each $ r \leq l \leq L$, $\left | \cJ_{z_{l+1}} \right | \geq \left | \cJ_{z_l} \right | + 1$ and we have $\left | \cJ_{z_L} \right | \leq n$. Remark \ref{prop:ILtoGIL_binding} shows the tradeoff in observation distance $r$ increases in $\GGIL$ leading to more constraints to bind. 

\begin{remark} \label{prop:ILtoGIL_binding}
Let $(\mathbf{v}_{r+1}, \bE_{r+1}, \bz_{r+1})$ be optimal for $\GGIL({\bX, \Omega, \bz_{r},\bP = \emptyset}),$ $\forall r \in \left \{ 1,\hdots,L \right \}$, then $\cD(\bE_{2}) \leq \cD(\bE_{3}) \leq  \hdots \leq \cD(\bE_{L})$.
\end{remark}

Despite the increase in the metric $\cD$ at optimality as the number of binding \desirable\ constraints grows, comparison between the solutions is easier since the sequence of solutions $\bz_1,\bz_{2}, \hdots, \bz_{L}$ are all contained in one face of $\Omega$.
%

\begin{remark} \label{prop:GGIL_Face}
Let $(\mathbf{v}_{r+1}, \bE_{r+1}, \bz_{r+1})$ be optimal for $\GGIL({\bX, \Omega, \bz_{r},\bP}),$ $\forall r \in \left \{ 1,\hdots,L \right \}$, then $\bz_{r} \in \mathcal{F} = \left \{ \bz \in \Omega | \ \ba_j \bz = b_j \ \forall j \in \cJ_{z_1} \right \},$ $\forall r \in \left \{ 1,\hdots,L \right \}$.
\end{remark}

Remark \ref{prop:GGIL_Face} shows that for a sequence of solutions from $\GGIL$, all are contained in a face of $\Omega$ by noting that Constraint \eqref{GGILBinding} ensures prior binding constraints will be tight at optimality. A subsequent result is that all the solutions in the sequence are optimal for any cost vector that makes the previous solutions optimal, as shown in Theorem \ref{Thm:GGIL_Optimal}. 

\begin{theorem} \label{Thm:GGIL_Optimal}
Let $\bz_1,\bz_{2}, \hdots, \bz_{L}$ be a sequence of learned solutions from successively solving $\GGIL$. Then, $\forall r \in \left \{ 1,\hdots,L \right \}$, $\bz_{r}$ is optimal for $\FO(\Bar{\bc}_l, \Omega)$ where $\Bar{\bc}_l \in cone (\ba_t : t \in \cJ_{z_l}), \ \forall l \in \left \{ 1,\hdots,r \right \}$ and $\Bar{\bc}_l \neq \bzero$.
\end{theorem}

Using the results of Theorem \ref{Thm:GGIL_Optimal}, we observe that by adaptively increasing the number of binding constraints in solutions using $\GGIL$, we achieve new solutions that are still optimal for the previously achieved cost vectors but bind additional \desirable\ constraints. 
This trait and the shared binding \desirable\ constraints among the solutions in a sequence provide easier comparison across the candidate learned solutions for the user. In other words, the user can assess the ``cost'' of binding more \desirable\ or \pref\ constraints.  

\subsection{Summary} \label{Sec:methods_summary}
The inverse learning methodology detailed in Sections \ref{Section:Methodology} and \ref{sec:AIL} is capable of incorporating existing knowledge of a linear optimization problem and known observations to learn optimal solutions, as well as recovering the missing parameters of the problem. 
The $\IO$ model incorporates information on constraints, if available, to learn optimal solutions that bind relevant constraints. 
The goal-integrated model, $\MGIL$, generalizes $\IO$ to control the \ObservationConstraintTradeoff\ through adjusting the number and the type of binding constraints. Finally, the $\GGIL$ models adapt a prior solution that is desirable to an improved new solution by preserving the previously binding constraints while adding to them. Iteratively solving  $\GGIL$ leads to a sequence of solutions that gradually increases achieved goals, hence, empowering users to balance and navigate the inherent tradeoff. 

We next validate our models and compare them with existing literature using a two-dimensional example problem in Section \ref{Section:NumericalEx}. In Section \ref{Section:Application}, we apply the $\GGIL$ model to the diet recommendation problem and demonstrate how comparable diet options can be derived for users. We show that by employing $\GGIL$ and creating a sequence of solutions starting with the solution that minimizes $\cD$, we can provide diverse diet recommendations with different yet desirable and measurable characteristics for the patients.

\section{Numerical Example} \label{Section:NumericalEx}

In this section, we consider an illustrative two-dimensional example to visually demonstrate the results of the inverse learning and the inverse optimization approaches. 
We consider a linear forward optimization problem with an unknown cost vector and seven known constraints, as outlined below, and a set of 20 feasible and infeasible observations (blue dots), as depicted in Figure \ref{Fig:exampleAllModels}. 
\\
$\FO(\bc,\Omega) = \underset{\bx\in\mathbb{R}^2}{max} \left \{ \bc'\bx~|~{\Gamma_1}: -x_1 + x_2 \leq 10, ~~~{\Gamma_2}: -0.5x_1 + x_2 \leq 11,\right.$\\ 
$\left.~~~~~~\hspace{2cm}{\Gamma_3}: 0.5x_1 + x_2 \leq 16, ~~~{\Gamma_4}: x_1 + x_2 \leq 20, ~~~{\Gamma_5}: x_1  \leq 10, ~~~{\Gamma_6}: x_1  \geq 0, ~~~{\Gamma_7}: x_2  \geq 0 \right \}$.
We assume constraints $\Gamma_1$--$\Gamma_5$ are \desirable\ with $\Gamma_3$ being a \pref\ constraint. We note that the example is designed in such a way that knowledge of \desirable\ constraints does not change the results of the benchmark literature models to facilitate a fair comparison. The \pref\ constraint is only added to showcase the goal-integrated inverse learning model ($\MGIL$) performance. For $\MGIL$, we consider three cases: binding one or two \desirable\ constraints (i.e., $r=1,2$) with and without prioritizing the \pref\ constraint, i.e., $\bP = \left \{3 \right \}$ and $\bP = \emptyset$, respectively.
%
We solve the same example problem with benchmarks inverse linear optimization models from the literature, namely, the absolute duality gap model \citep{babier2021ensemble, chan2019inverse} and the mixed-integer quantile inverse optimization model \citep{shahmoradi2021quantile}, also referred to as $\ADG$  and $\MIIO$ for brevity, respectively. 
For the duality gap model, we solve the recovered $\FO(\bc,\Omega)$ using the Simplex method to find an optimal solution. For the quantile model, we consider the quantile error to be $\theta =1$ as suggested by the authors alongside depicting results for $\theta =0.9$ and $\theta =0.8$. To find the best possible objective value, we devised a bisection algorithm to find the minimum threshold optimality error $\tau$, as no default value or suggested method is provided in the work. 
For all models and whenever appropriate, we used the 2-norm distance metric to learn $\FO$ optimal solutions and to recover the best cost vector. 

\iftrue
\begin{figure}[htb]
    \begin{center}
        \includegraphics[width=.6\linewidth]{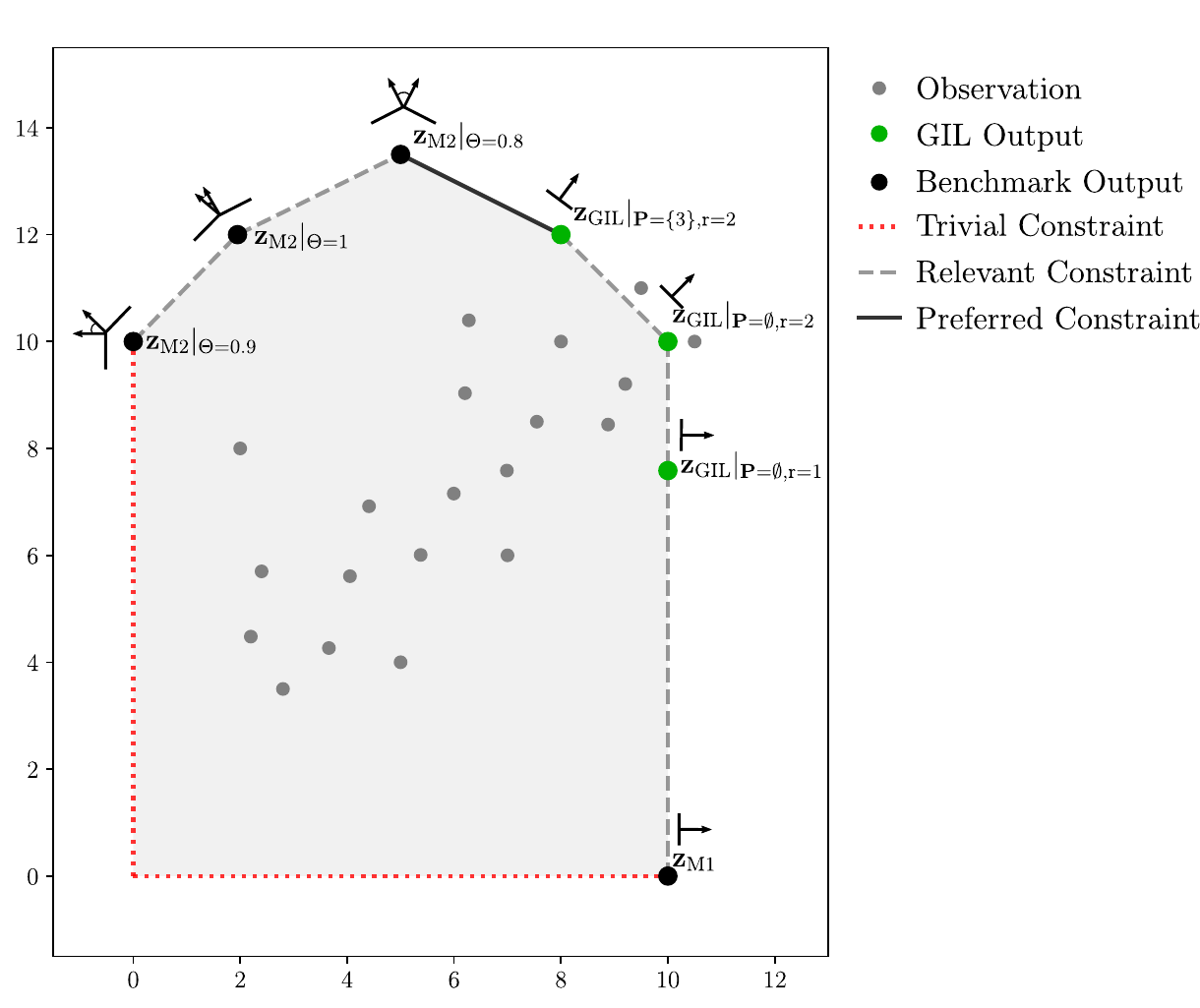}
    \caption{An illustrative 2D example with 20 observations. The recovered optimal solutions and cost vectors of $\FO$ are shown for the absolute duality gap model ($\ADG$), the quantile model ($\MIIO$), and the goal-integrated inverse learning model ($\MGIL$). The observations demonstrate a trait of recurring closer to the top-right constraint bounds. $\MGIL$ captures this trend while the solutions from $\ADG$ and $\MIIO$ models can be far from the observations. If additional constraint knowledge is included, $\MGIL$ can incorporate them to adjust the learned solution and cost vector accordingly.}
    \label{Fig:exampleAllModels}
    \end{center}
\end{figure}
\fi

\iftrue
\begin{table}[h]
\footnotesize
\begin{center}\renewcommand{\arraystretch}{1.2}
        \caption{\footnotesize  Results from applying the inverse learning model and the existing models in the literature: absolute duality gap ($\ADG$) and the quantile models ($\MIIO$). The solutions learned from inverse learning models exhibit reduced distances to the observations and recover similar cost vectors.}
    \label{Table:DietCopmILvsGIL}
    \begin{tabular}{>{\centering}p{0.17\textwidth}|>{\centering}p{0.15\textwidth}|>{\centering}p{0.16\textwidth}|>{\centering}p{0.11\textwidth}|>{\centering}p{0.15\textwidth}|>{\centering\arraybackslash}p{0.14\textwidth}}
    \hline\hline
    \footnotesize
    Model  &  \makecell{Avg. Distance \\to Observations}  & \makecell{Set of $\FO$ \\Cost Vectors}  &
    \makecell{$\FO$ \\Cost Vector}& \makecell{$\FO$ \\Optimal Solution}  & \makecell{Avg. Objective \\Value$^*$} \\
    \hline 
    \multicolumn{6}{c}{\textbf{Inverse Learning Models}}\\
    \hline
    $\MGIL(r = 1, \bP = \emptyset)$ & 6.0 &  C**((1,0)) & (1,0) & (10, 7.58) &  5.9 \\
    $\MGIL(r = 2, \bP = \emptyset)$ & 7.0 &  C((1,0),(1,1)) & (0.7,0.7) &   (10, 10) &  9.2 \\
    $\MGIL(r=2,\bP = \left \{3 \right \})$ & 7.4 &  C((1,1),(0.5,1)) & (0.63,0.78) &  (8, 12) &  9.4 \\
    \hline
    \multicolumn{6}{c}{\textbf{Benchmark Inverse Optimization Models}}\\
    \hline
    $\ADG$ & 11.4 &  (1,0) & (1,0) &  (10, 0) &  5.9 \\
    $\MIIO{(\theta=1)}$ & 8.7 &  C((-1,0),(-0.5,1)) & NA$^{***}$ &  (1.95, 12) &  NA \\
    $\MIIO{(\theta=0.9)}$ & 8.7 &  C((-1,0),(-1,1)) & NA &  (0, 10) &  NA \\
    $\MIIO{(\theta=0.8)}$ & 8.5 &  C((-0.5,1),(0.5,1)) & NA &  (5, 13.5) &  NA \\
    \hline\hline
    \end{tabular}
\end{center}
*Average forward objective values for all the observations using the recovered cost vector in the forward problem. \\ 
**Represents the cone of vectors. \\
*** Not applicable. The reference study does not provide a direct method to select a single cost vector.
\end{table}
\fi

As Figure \ref{Fig:exampleAllModels} illustrates, the recovered solution by the absolute duality gap model ($\ADG$), as an example of classical inverse optimization, may be far from the observations and unstable, as known in the literature. 
The solutions found by the quantile model ($\MIIO$), while more stable, can also be located far from the observations and are sensitive to the parameters $\theta$ quantile and the minimum threshold optimality error $\tau$. If perfect values for these parameters are known, it is possible to guide the solution closer to the observations but finding such values is not trivial and may also not be intuitive to the end-users to express them in the required form.   
The learned solutions of $\MGIL$, however, can resemble the observations' trends. The optimal solution returned by $\MGIL(\bX,\Omega,r=1,\bP = \emptyset)$ (which is equivalent to the $\IO$ model) is the closest to the observations. The location of the learned solution can be further controlled by requiring two relevant constraints to be active, as shown in $\MGIL(\bX,\Omega,r=2,\bP = \emptyset)$ which returns the closest extreme point to the observations. We can further tailor the solution by prioritizing constraint $\Gamma_3$ to bind as our \pref\ constraint by setting $\MGIL(\bX,\Omega,r,\bP = \left \{3 \right \} )$ which returns the closest feasible point on $\Gamma_3$ but at a potential cost of moving away from the observations.  
While it is possible to characterize the set of all inverse optimal cost vectors using $\MGIL$ models (as demonstrated in Theorem \ref{Theorem2}), using formulation \eqref{InverseOptimalc} allows $\MGIL$ models to find cost vectors that do not solely depend on $\Omega$, and that are informed of the observations as well. 

Table \ref{Table:DietCopmILvsGIL} summarizes the results shown in Figure \ref{Fig:exampleAllModels}. The set of returned optimal cost vectors, a cone except for $\ADG$ that is a point, the selected cost vector from said cones, for frameworks that identify it, and the optimal solutions are provided. The average distance to the observations is the smallest in $\GIL$ although it can increase as the solutions are tailored to specific constraint goals. 
The last column of Table \ref{Table:DietCopmILvsGIL} indicates how using formulation \eqref{InverseOptimalc} results in maximizing the average objective value of the observations for $\FO$. 

\section{Application to Personalized Diet Recommendations} \label{Section:Application}

In this section, we apply the inverse learning models to a personalized diet recommendations problem using observations from 2,590 hypertension patients and develop a web-based decision-support tool for exploring the data and the results. The diet recommendation problem is among the most well-studied problems in operations research~\citep{stigler1945cost,bassi1976diet,dantzig1965linear,gazan2018mathematical}. The goal of the classical diet recommendation problem, as first introduced by \cite{stigler1945cost}, is to find the optimal food intake of a patient or a set of patients based on some given cost function subject to a given set of constraints \citep{stigler1945cost,garille2001stigler}. Many studies have focused on finding optimal diets for individuals and communities based on existing behaviors  \citep{stigler1945cost, sacks2001effects,bas2014robust}. The more recent literature has focused on providing more realistic constraints and relaxing previous simplifying assumptions of the original diet recommendation problem, such as the linear relationship between the nutrients and foods and the varying effects of different combinations of foods \citep{bas2014robust}.

In the inverse optimization literature, the diet recommendation problem has been studied to recover missing cost functions (or equivalently, the cost functions)~\citep{ghobadi2018robust,shahmoradi2021quantile} or constraint parameters \citep{ghobadi2021inferring} based on a given set of observed dietary decisions. In our case study, we concentrate on recommending diets and recovering a corresponding cost function based on a given set of observed daily dietary decisions. Similar to \cite{ghobadi2018robust}, in the definition of the diet recommendation problem, we assume that the cost function is inversely associated with the preferences of the individual, meaning that the more palatable a particular food item for an individual, the smaller its cost. 

In this work, we consider the dietary behaviors of individuals who self-reported hypertension. 
Roughly one out of every two adults in the United States experiences some level of hypertension or high blood pressure \citep{CDC_hypertension_2020}. For such patients, many clinicians recommend the Dietary Approaches to Stop Hypertension (DASH) eating plan that aims to decrease sodium intake and control blood pressure \citep{challa2018dash}. The DASH diet has been shown to lower both systolic and diastolic blood pressures in individuals \citep{sacks2001effects,challa2018dash}. 
However, maintaining a long-term and realistic diet is a persistent challenge for many patients. Diets with clear options that resemble the existing dietary behaviors of the patients and have increased palatability are more likely to be adhered to by the patient for long-term commitments \citep{bentley2005factors}. Many of the existing methods in the literature rely on defining specific constraints in hopes of achieving more palatable diets \citep{ferguson2006design,xu2018cs}. Here, we consider a forward optimization model based on the DASH diet nutritional constraints and the palatability objective. We then utilize the inverse learning models using given dietary behavior data of patients to recommend diets ranging from diets similar to the decisions of the patients while adhering to DASH guidelines to quality diets that are richer in many nutrients. 

In this section, we first discuss the details of the data sets used for our application of the diet recommendation problem in Section \ref{sec:app_data}. We then showcase the results of applying inverse learning for the purpose of recommending group-preferred diets considering nutritional limits in Section \ref{sec:app_results}. Finally, we introduce a web-based online tool in Section \ref{sec:webpage} that can be used to provide and compare inverse learning recommendations for different groups of patients in our dataset.

\subsection{Data} \label{sec:app_data}
The input to inverse learning models includes a set of observations and a set of constraints that form the known feasible set $\Omega$. We use databases from the National Health and Nutrition Examination Survey (NHANES) \citep{CDC_2020} and the United States Department of Agriculture (USDA) database \citep{usda_2019} to generate the observations and the feasible set of the forward optimization problem, respectively. NHANES database includes self-reported information on the daily food intakes of 9,544 individuals, along with their demographics, for two days. From this population, 2,090 individuals self-reported a diagnosis of hypertension totaling to 4,024 observations of daily food intake for patients who did report their daily food intakes for two days.  
We consider the DASH diet as the recommended diet for our $\FO$ problem and build lower-bound and upper-bound nutritional constraints based on this diet recommendation for patient clusters based on their age and gender ~\citep{liese2009adherence,sacks2001effects}. 
To calculate the nutrients per serving of each food group and construct the feasibility set, we employ information from the USDA databases. The USDA datasets are detailed with approximately 5,000 food types listed. For simplicity and tractability in our models, we group the food types into 38 food groups as shown in Table \ref{Table:food_groups} in the electronic companion. We additionally consider 22 nutritional constraints and 76 box constraints on food types and limit the intake amount for each food item to a maximum of eight servings per day.

\iftrue
\begin{table}[htp]
\small
\begin{center}\renewcommand{\arraystretch}{1.2}
        \caption{DASH nutrient bounds for women of 51+ age group. The nutritional coefficients per serving of foods for a sample of food groups are also provided.}
    \label{Table:SubsetNutrients}
    \begin{tabular}{>{\centering}p{0.19\textwidth}|>{\centering}p{0.1\textwidth}|>{\centering}p{0.1\textwidth}|>{\centering}p{0.08\textwidth}>{\centering}p{0.08\textwidth}>{\centering}p{0.08\textwidth}>{\centering\arraybackslash}p{0.19\textwidth}}
& \multicolumn{2}{c|}{\makecell{ \textbf{DASH Nutritional} \\ \textbf{Bounds}}}  &  \multicolumn{4}{c}{\makecell{ \textbf{Nutrient Per Serving for Select Foods} }} \\\cline{2-7}
\textbf{Nutrients}        & \textbf{Lower bound} & \textbf{Upper bound} & \textbf{Milk (244g)}   & \textbf{Stew (140g)}   & \textbf{Bread (25g)} & \textbf{Tropical Fruits** (182g)}  \\
\hline \hline
Energy (kcal)    & 1,575.1$^*$     & 2,013.2$^*$     & 160.1 & 282.0 & 72.5 & 120.1 \\
Carbohydrate (g) & 223.7$^*$      & 254.2      & 18.8  & 21.3  & 12.3 & 27.8     \\
Protein (g)      & 51.8       & 89.2$^*$      & 7.2   & 16.8  & 2.5 & 1.4      \\
Total Fat (g)    & 59.7$^*$       &  78.2      & 6.3   & 14.3  & 1.5 & 1.9      \\
Total Sugars (g) & 117.0$^*$       &  144.6      & 18.0  & 4.5   & 1.7 & 19.0     \\
Fiber (g)        & 36.7       &  39.3$^*$       & 0.2   & 1.5   & 0.9 & 4.1      \\
Sat. Fat (mg)    & 11.4$^*$        & 16.6       & 3.3   & 4.6   & 0.4 & 0.31      \\
Cholesterol (mg) & 24.4$^*$       &  120.6      & 14.0  & 53.9  & 1.2 & 0.0      \\
Iron (mg)        & 9.7        & 12.5$^*$       & 1.0   & 2.2   & 0.8 & 0.5      \\
Sodium (mg)      & 1,376.2$^*$     &  1,693.0     & 108.3  & 639.4 & 119.0 & 10.7    \\
Caffeine (mg)    & 0.0$^*$        & 80.0$^*$       & 0.3   & 0.0   & 0.0 & 0.0    \\
    \hline  \hline
    \end{tabular}
\end{center}
* Indicates \desirable\ constraint. \\
** Includes apples, apricots, avocados, bananas, cantaloupes, cherries, grapes, mangoes, and pineapples.
\end{table}
\fi

\iftrue
\begin{table}[]
\small
    \caption{\footnotesize Description of the nutritional facts of the observations for patients alongside minimum and maximum value of each nutrient among observations with 25, 50, 75$^{\text th}$ percentiles.}
    \label{Table:Diet_observations}
\begin{center}\renewcommand{\arraystretch}{1.2}
\begin{tabular} {>{\centering}p{0.19\textwidth}|>{\centering}p{0.1\textwidth}|>{\centering}p{0.1\textwidth}|>{\centering}p{0.08\textwidth}>{\centering}p{0.08\textwidth}>{\centering}p{0.08\textwidth}>{\centering}p{0.08\textwidth}>{\centering\arraybackslash}p{0.08\textwidth}}
\hline\hline
\textbf{Nutreint}  &  \textbf{Avg.} & \textbf{Std} & \textbf{Min} & \textbf{25\%} & \textbf{50\%} & \textbf{75\%} & \textbf{Max} \\ \hline \hline
Food Energy (kcal) &  1,889.8        & 787.1        & 400.5        & 1,329.3       & 1,780.9       & 2,323.6       & 4,969.5       \\
Carbohydrate (gm)  &  242.5         & 105.3        & 48.1         & 167.7        & 231.6        & 299.2        & 647.9        \\
Protein (gm)       &  69.6          & 34.1         & 10.6         & 46.4         & 64.1         & 85.3         & 233.9        \\
Total Fat (gm)     &  74.6          & 36.8         & 6.4          & 48.7         & 69.3         & 93.2         & 239.4        \\
Sugars (gm)        &  104.9         & 54.8         & 16.3         & 65.0         & 96.9         & 131.9        & 446.6        \\
Dietary Fiber (gm) &  15.6          & 8.8          & 1.5          & 9.7          & 13.6         & 19.8         & 54.9         \\
Saturated Fat (gm) &  23.0          & 11.5         & 1.6          & 14.6         & 21.2         & 29.9         & 72.5         \\
Cholesterol (mg)   &  250.2         & 192.9        & 5.5          & 124.8        & 196.9        & 331.7        & 1,262.0       \\
Iron (mg)          &  13.6          & 7.7          & 2.3          & 8.3          & 12.1         & 16.5         & 52.9         \\
Sodium (mg)        &  3,413.4        & 1,643.2       & 324.7        & 2,205.5       & 3,106.6       & 4,288.4       & 9,942.8       \\
Caffeine (mg)      &  138.3         & 119.1        & 0.0          & 54.9         & 114.3        & 190.4        & 1,021.3    \\   \hline \hline
\end{tabular}
\end{center}
\end{table}
\fi

The complete dataset of the self-reported dietary decisions of patients in NHANES is highly heterogeneous and the DASH diet recommends different target calories based on age and gender. 
As a result, we consider a simple demographic clustering based on the age and gender of the participants. For illustrative purposes, we focus on a demographic group of women who are 51+ years old, are self-reported both hypertension and pre-diabetes, and have at least two observations per patient. These restrictions will reduce our sample population to 230 patients and 460 observations. Other demographic groups can be explored in our web-based tool described in Section~\ref{sec:webpage}.

Table~\ref{Table:SubsetNutrients} presents the DASH diet bounds on nutrients for our representative demographic group. The table also indicates the \desirable\ constraints and demonstrates the coefficients for a sample of food items. All nutrients include at least one desirable bound which is considered a \desirable\ constraint. For the cases where we do not assume any preference for the nutrient, both bounds are considered as \desirable\ constraints. 
Table~\ref{Table:Diet_observations} summarizes the nutritional values of the 460 observed dietary decisions corresponding to our representative patient group (since each patient reports two days of food intake). A quick comparison with Table~\ref{Table:SubsetNutrients} shows that a large number of observations are infeasible, especially for the sodium constraint which the DASH diet aims to restrict. For instance, in the representative group of 51+ year-old women, roughly 80\% of observations are infeasible with sodium causing the most number of infeasibilities (~70\%). 
These infeasible observations still contain information about patient preferences, and hence, we include them in our inverse learning models. 
For further information about the constraints, food types, observed dietary decisions, and other relevant information, please see Section \ref{Appendix:diet_data} in the electronic companion. 

\subsection{Recommending Diets with Inverse Learning} \label{sec:app_results}
In this section, we demonstrate the use of inverse learning models for diet recommendations and use the adaptive $\GGIL$ models to enable navigating between different diet recommendations. 
The decision-maker can inform the models of their \pref\ nutritional constraints, if any, and learn recommended diets that balance between distance to prior observed dietary decisions and binding \pref\ constraints. 
For comparison purposes, we consider two settings of having no \pref\ constraints and having explicit \pref\ constraints. 
As an example of \pref\ constraints, we consider four constraints as \pref\ by the decision-maker: the lower bounds of sodium, saturated fat, and cholesterol, and the upper bound of dietary fiber. We denote the vector of their corresponding indices by $\Bar{\bP}$. A user can choose other constraints as \pref\ constraints using the web-based tool that is described in Section \ref{sec:webpage}.

We first solve $\MGIL(\bX,\Omega,r=1, \bP = \emptyset)$ and $\MGIL(\bX,\Omega,r=1, \bP = \Bar{\bP})$, which similar to $\IL$ models, bind one (or more) relevant constraints and denote the learned diets as $\bd^{\emptyset}_1$ and $\bd^{\Bar{\bP}}_1$, respectively. 
The subscripts indicate the minimum number of binding \desirable\ constraints ($r$). 
We then generate two sequences of adaptive diets by iteratively solving the models for $r=1,\dots, 7$, denoted by $\bd^{\emptyset}_1,\hdots,\bd^{\emptyset}_7$ and $\bd^{\Bar{\bP}}_1,\hdots,\bd^{\Bar{\bP}}_7$, respectively. 
For all models, we consider the 1-norm distance metric and the demographic group of 51+ year-old women. 
Other demographics can be explored on the web-based tool described in Section \ref{sec:webpage}. 

\begin{figure}[t]
\begin{center}
\includegraphics[width =.8 \linewidth]{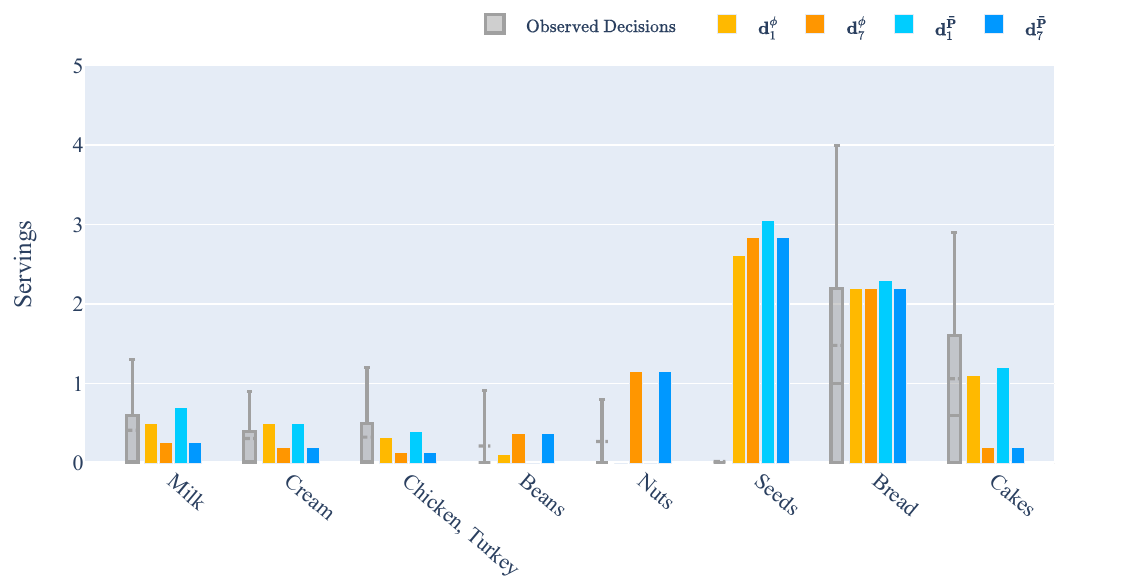}
\includegraphics[width =.8 \linewidth]{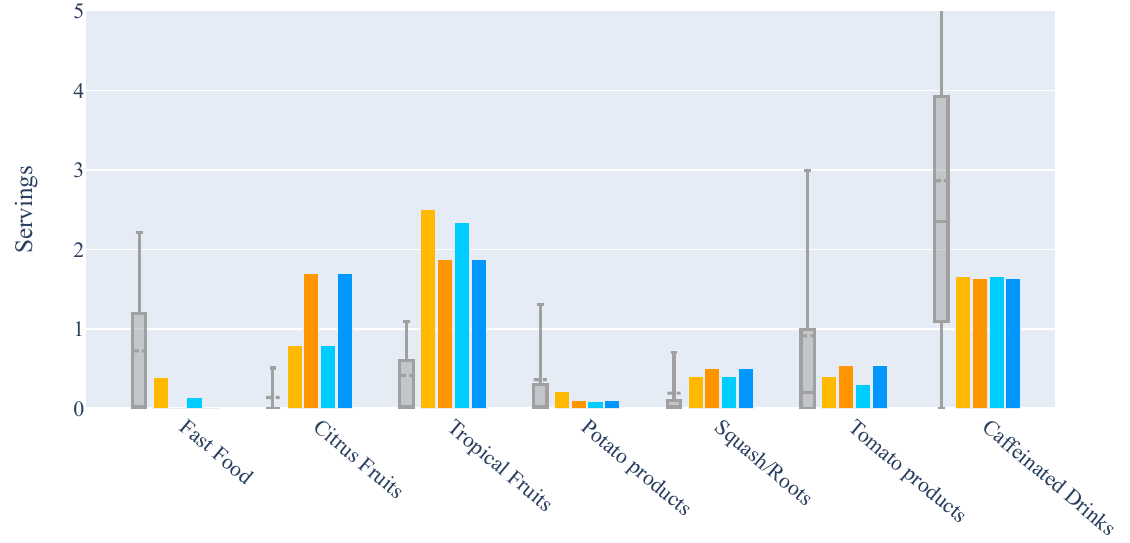}
\caption{\footnotesize A set of 460 historical observations of daily diets for the 51+ year-old women (boxplots) is shown with a set of four learned diets with binding constraints ($r=1$ or $7$) and with or without preferred constraints ($\bP=\emptyset$ or $\bar{\bP}$). The learned diets tend to preserve historical behavior when $r=1$ and move further from observations and towards higher food value identified by DASH when $r$ increases or a set of preferred constraints are identified. 
} \label{Fig:ILvsMILfoods}
\end{center}
\end{figure}

Figure \ref{Fig:ILvsMILfoods} illustrates four representative learned diets of  $\bd^{\emptyset}_1$ $\bd^{\emptyset}_7$, $\bd^{\Bar{\bP}}_1$, and $\bd^{\Bar{\bP}}_7$, which bind $r=1$ or $7$ constraints with or without the specified preferred constraints  $\bar{\bP}$ (fiber, sodium, cholesterol, and saturated fat). 
The historical daily intakes, which serve as observations in our models, are shown in gray boxplots (at 10, 25, 50, 75, and 90$^{\text th}$ percentiles) for a sample set of food items (full list in Section \ref{Appendix:diet_data}). 
The prior behavior of this demographic group is mostly infeasible for DASH requirements in fiber and sodium levels and shows, for instance, small levels of consumption of nuts and seeds but higher levels for bread and caffeinated drinks. 
The learned diets trend toward increased values in healthier food items (including nuts and fruits) due to the DASH restrictions. 
For example, the learned diets recommend an increased intake of nuts and seeds (unlike the observations), a reduced intake of fast food (also unlike the observations), but similar levels of consumption for bread, milk \& cream to the observations. 
We note that the caffeinated drinks are capped in all learned diets due to an upper bound on caffeine in DASH.

\begin{figure}[t]
\begin{center}
\includegraphics[width =.8 \linewidth]{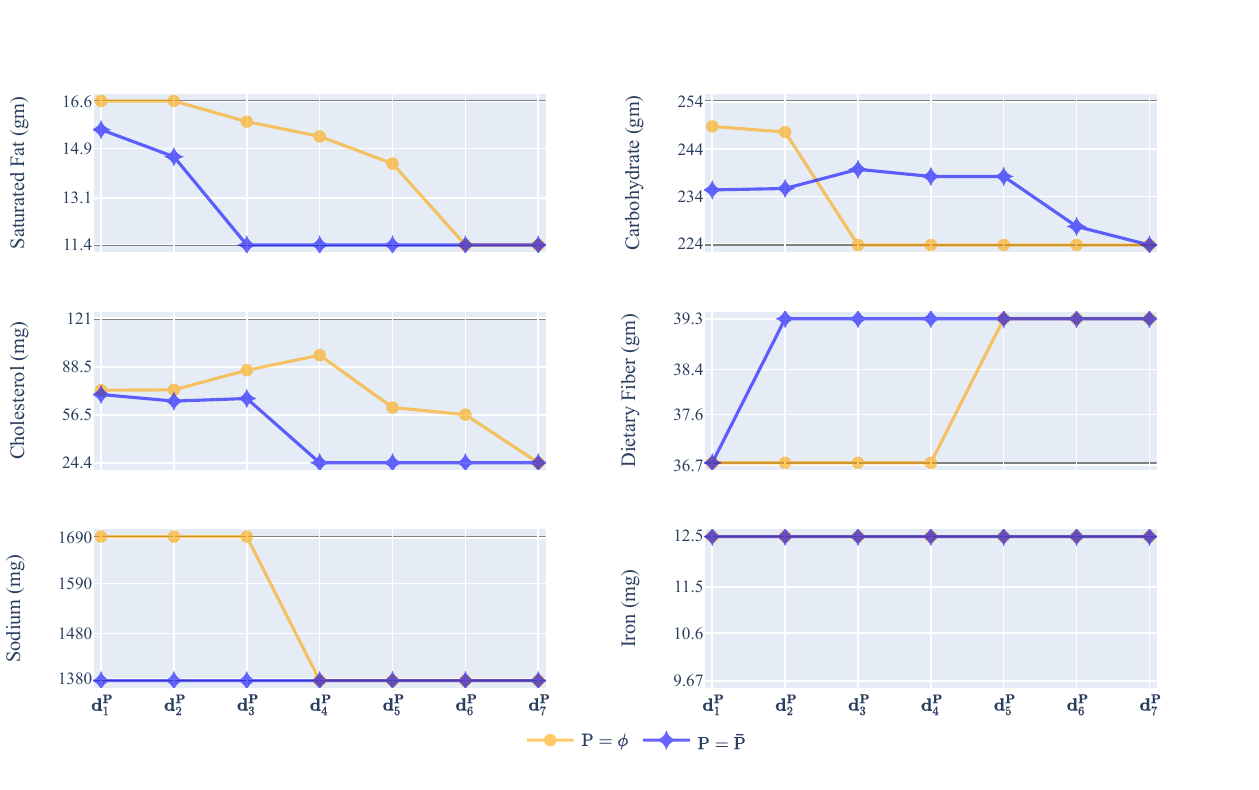}
\caption{\footnotesize A select set of DASH nutritional constraints for two adaptive sequences of $\bd^{\emptyset}_r$ (yellow) and $\bd^{\Bar{\bP}}_r$ (blue), $r=1,\dots,7$, where $\bar{\bP}=$ Saturated Fat, Cholesterol, Sodium, and Dietary Fiber. 
The \pref\ constraints bind faster when they are specified explicitly ($\bP = \bar{\bP}$). For relevant constraints that are not part of the \pref\ constraints, e.g., carbohydrate, they may or may not bind depending on the given value of $r$. 
} \label{Fig:ILvsMILnuts}
\end{center}
\end{figure}

The learned diets $\bd^{\emptyset}_1$ and $\bd^{\Bar{\bP}}_1$, which bind one relevant (or preferred) constraint, consistently preserve the historical behavior of the demographic group better (closer to the boxplots) than when more constraints are binding, for instance, see intakes for cakes. Increasing the number of binding constraints or specifying more \pref\ constraints can increase the distance from the observations (as also illustrated in Table \ref{Table:ILvsMIL}). 
Specifically, diets recommended by considering the non-empty set of \pref\ constraints $\Bar{\bP}$ show a sharper trend towards food groups richer in dietary fiber and food groups containing less sodium, the two main criteria in DASH that most historical observations of this demographic group do not satisfy. 

Similar trends can also be seen for nutrient constraints as shown in Figure \ref{Fig:ILvsMILnuts}. The figure compares the amounts of select nutrients for the learned diets in the two sequences of $\bd^{\emptyset}_r$ and $\bd^{\Bar{\bP}}_r$, $r=1,\dots,7$. 
Inverse learning models recommend healthier diets with reduced amounts of fats for the case where $\bP = \emptyset$ and reduced amounts of carbohydrates, fats, cholesterol, and sodium for the case where $\bP = \Bar{\bP}$. 
The sequence of solutions obtained for $\Bar{\bP}$ binds \pref\ constraints faster (at upper- or lower-bound, whichever that is specified) than the case where we do not consider \pref\ constraints. This trend showcases how considering \pref\ constraints prioritize healthier diets over replicating prior dietary choices of the patients. When a relevant constraint is not part of the \pref\ constraints, it may bind for all models (e.g., upper bound of iron) or bind in $\IL$ faster than $\GIL$ (e.g., lower bound of carbohydrate). 
Expanded results for all nutrients are included in the electronic companion. 

\iftrue
\begin{table}[]\renewcommand{\arraystretch}{1.2}
\small
\footnotesize
    \centering
        \caption{\footnotesize Two learned adaptive sequences of diets. The Average distance to observation is non-decreasing for increased values of the parameter $r$ while the number of binding \desirable\ and \pref\ constraints increases.}
    \label{Table:ILvsMIL}
    \begin{tabular}{>{\centering}p{0.085\textwidth}|>{\centering}p{0.08\textwidth}>{\centering}p{0.08\textwidth}|>{\centering}p{0.1\textwidth}>{\centering}p{0.08\textwidth}|>{\centering}p{0.08\textwidth}>{\centering}p{0.08\textwidth}|>{\centering}p{0.08\textwidth}>{\centering\arraybackslash}p{0.08\textwidth}}
    \hline\hline
\centering
\multirow{2}{5em}{Learned Solution} & \multicolumn{2}{c}{\makecell{ Avg. Distance \\to Observations}} & \multicolumn{2}{c}{\makecell{\# of Binding \\\Desirable\ Constraints}}     & \multicolumn{2}{c}{\makecell{\#  of Binding \\\Pref\ Constraints}} & \multicolumn{2}{c}{\makecell{\#  of Binding \\\Objective\ Constraints}}\\ 
                              & $\bP = \emptyset$ & $\bP = \Bar{\bP}$ & $\bP = \emptyset$ & $\bP = \Bar{\bP}$ & $\bP = \emptyset$ & $\bP = \Bar{\bP}$ & $\bP = \emptyset$ & $\bP = \Bar{\bP}$ \\ \hline
                           $\bd^\bP_1$  & 18.4  & 18.7  & 1  & 3  & 0 & 1 & 6 & 4  \\
                           $\bd^\bP_2$  & 18.4  & 18.9  & 2  & 4  & 0 & 2 & 6 & 3 \\
                           $\bd^\bP_3$  & 18.5  & 19.2  & 3  & 5  & 0 & 3 & 5 & 2 \\
                           $\bd^\bP_4$  & 18.7  & 19.7  & 4  & 5  & 1 & 4 & 4 & 2 \\
                           $\bd^\bP_5$  & 19.0  & 19.7  & 5  & 5  & 2 & 4 & 3 & 2 \\
                           $\bd^\bP_6$  & 19.2  & 19.7  & 6  & 6  & 3 & 4 & 2 & 2 \\
                           $\bd^\bP_7$  & 19.7  & 19.7  & 7  & 7  & 4 & 4 & 2 & 2 \\
    \hline \hline
    \end{tabular}
\end{table}
\fi

Table \ref{Table:ILvsMIL} summarizes the results for the two adaptive sequences of recommended diets. The average 1-norm distance for all recommended diets to the historical observations is presented, as well as the number of binding nutritional \desirable, \pref, and \objective\ constraints. As expected, the average distance of the learned diet increases as the number of binding constraints $r$ increases as the learned diet may be further from the observed dietary behaviors as well as the number of binding \desirable\ and \pref\ constraints. The number of binding \objective\ constraints, however, decreases in both sequences with a faster decrease when $\bP=\bar{\bP}$. 
We note that over 80\% of all historical observations in this example are infeasible for DASH requirements, and hence, the initial distance for $r=1$ is large. 
In general, including \pref\ constraints increases the distance to the observations compared to the cases where $\bP = \emptyset$. However, in the case of this particular example, we observe that the increase in distance from the observations is not drastically higher. This tradeoff can be a measure for the user to opt for recommending the sequence of solutions obtained from solving $\GGIL(\bX,\Omega,\bd^{\Bar{\bP}}_r,\bP = \Bar{\bP})$ to the patient as they bind more \pref\ constraints without increasing the distance by a large margin. It rests with the user/decision-maker to indicate which of the two competing parameters is closer to their dietary goals. 
Finally, although computational efficiency is not a focus of this work, we note that the computational time for inverse learning models is comparable with existing models in the literature.

\subsection{Diet Recommendation Decision-support Tool} \label{sec:webpage}

The inverse learning framework enables the integration of additional knowledge (provided by domain experts or derived from data). Exploring a set of results allows the user to select the best solution or a sequence of solutions that best fit their needs. Given this flexibility and range, an interactive and user-friendly decision-support tool can empower users further. To this end, we develop a web-based interactive dashboard for our case study of diet recommendations and the National Health and Nutrition Examination Survey data \citep{CDC_2020}. The tool is available publicly at \href{https://optimal-lab.com/optimal-diet}{https://optimal-lab.com/optimal-diet} and includes all data and results presented in Section \ref{Section:Application}. 
Users can select different demographic groups, choose their own preferred constraints, and identify the number of binding constraints to investigate results and contrast them against each other. 

Figure \ref{Fig:website_IL} shows a screenshot of the interactive tool. The demographic group of choice can be selected from the dropdown menu marked as (a), for example, males between 31-50 years of age. The resulting diet recommendation for this group is then presented in part (b) with the binding lower- or upper-bound constraints displayed on the right, e.g., for this example sugar and iron. The historical observations are shown in boxplots with the average distance to observations shown on top left of the resulting figure. The diet recommendations are obtained by solving the $\GGIL$ models, and the user can explore binding different numbers of constraints by using the slider in part (c). This example illustrates the results for binding two constraints ($r=2$). The user can then define \pref\ constraints to bind from a set of DASH requirements as shown in part (d) and rerun the models with their new selection. 

\iftrue
\begin{figure}[tb]
\begin{center}
\includegraphics[width =0.95 \linewidth]{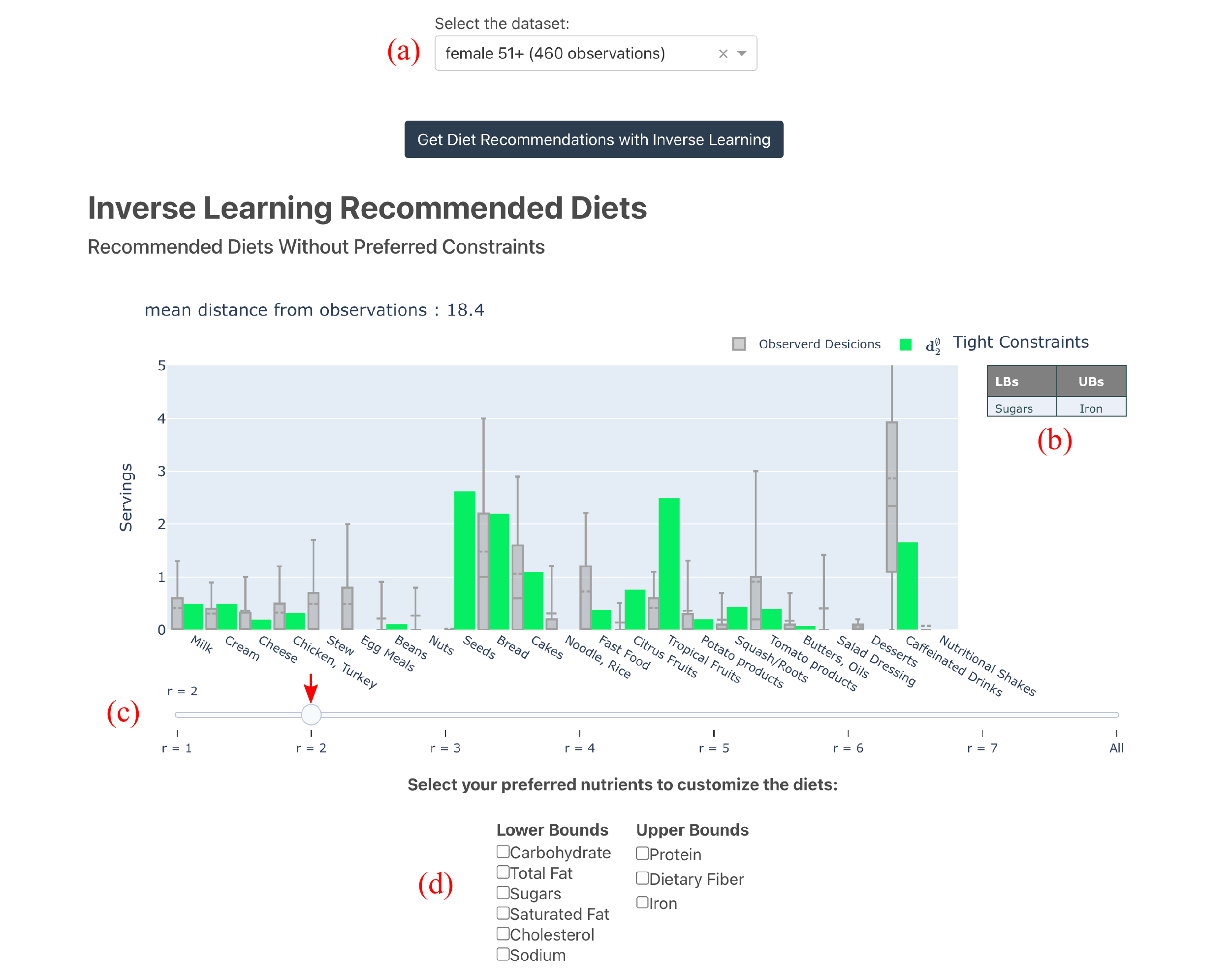}
\caption{\footnotesize A snapshot of the interactive decision-support tool that enables exploring a range of possible diets. 
Part (a) shows the dataset selection, (b) displays the binding relevant constraints, (c) enables exploring different numbers of binding constraints, and (d) provides a panel to choose \pref\ constraints. 
} \label{Fig:website_IL}
\end{center}
\end{figure}
\fi

The tool illustrates each learned diet with the average observed food intake. For each new model, a new set of figures are produced that similar to Figures \ref{Fig:ILvsMILfoods} and \ref{Fig:ILvsMILnuts}, showcase the details of the recommended diet with additional information on its nutritional bounds. 
The nutritional bound graph consists of eleven sub-graphs, each representing one of the nutrients considered in the model. The y-axis on each sub-graph is scaled to the range of lower and upper bounds on the nutrient. The range limits are indicated by a dashed green line for the \desirable\ constraints, a blue curve for the \pref\ constraints, and a yellow curve for other constraints. 
Should the user decide to run a second model with a different selection, they can see their new results alongside the previous one for easier comparison. The results for the latest two selections can be navigated using  \textit{Current Results} and \textit{Previous Results} sections.


The interactive decision-support tool enables users to access a range of recommended diets and explore different settings and preferences. They will be able to examine the impact of binding more goals (achieving more nutritional goals) in a step-by-step manner and evaluate the gains and the price, in the form of moving away from their habits and prior selections. The online tool contrasts with one-size-fits-all dietary approaches. It expands diet recommendations from a single fixed diet to a flexible range of diets that are informed by the user, which may improve adherence. 

\section{Conclusions} \label{Section:Conclusion}

In this work, we introduce the inverse learning framework as a data-driven inverse optimization method that learns optimal solutions to linear optimization problems with unknown cost vectors given a set of observed solutions. Inverse learning can be used in a wide range of applied settings based on historical observations or decisions. 
This new framework unifies the recovery of forward cost vectors and optimal solutions together and provides flexibility in learning optimal solutions by considering both feasible and infeasible observations. The models integrate additional knowledge on the importance of constraints at optimality. The generalized goal-integrated inverse learning models can tailor the optimal solutions based on the number and the type of active constraints at optimality and prioritize binding those that are preferred by the user or based on the problem setting. Adaptive models are developed to preserve properties of prior solutions while guiding them towards more desired ones.  
The framework provides a range of solutions that are interpretable and enables users to navigate a range of solutions to find the most suitable one for their application. The models provide a mechanism to balance the inherent tradeoff between retaining characteristics of the historical observation and binding specific \desirable\ constraints for a more contoured solution. 
We apply the inverse learning framework to derive diet recommendations based on historical dietary behaviors of hypertension patients. 
The models recommend a flexible set of diets based on observations and can tailor to potential preferences of the users for nutritional goals or desirable food groups. The recommended diets honor the healthy choices of the patients and replicate their favored dietary choices, when nutritional constraints permit, and provide a range of options for improving the behavior based on the user's priorities.  
%
We also showcase the results of the inverse learning models alongside benchmark methods on a two-dimensional representative numerical example. 
The results exhibit increased control over the distance to the observations and decreased sensitivity to potential outliers as perturbations of all observations are taken into account in the objective of the inverse learning models.

Given the data-driven nature of the inverse learning framework, the quality of data can impact the results, and it is recommended to consider leveraging statistical or data science approaches in combination with inverse learning. Developing data-handling methods and safeguarding against uncertainties, errors, biases, and outliers in the data can be a future direction.   
Additionally, the methods provided here are built on the assumption that the feasible set of $\FO$ is known. If in addition to the unknown cost vector, part of the feasible region is also unknown, new methods can be developed to recover both parameters. 


\OneAndAHalfSpacedXI
\bibliographystyle{informs2014} 
\bibliography{IO_multipoint} 
\clearpage

\ECSwitch
\ECHead{Electronic Companion}

\section{Dual-Form of the Forward Optimization Problem} \label{Sec:Dual_form}
Considering that inverse models are usually defined based on optimality constraints of the forward optimization problem using duality theorems, we include the following dual form of $\FO$ here as follows:
\begin{subequations} \label{DFO}
\begin{align}
\DFO: \underset{\by, \bu}{\text{maximize}} & \quad \by' \bb + \bu' \bh, \label{eq:DFO_obj}\\
\text{subject to} 
& \quad \bA' \by + \bG' \bu  =\bc, \label{DFODualFeas1}\\ 
& \quad \by \geq \bzero, \quad \bu \geq \bzero. 
\end{align}
\end{subequations}
In the above formulation, $\by$ is the vector of dual variables associated with the \desirable\ constraints and $\bu$ is the vector of dual variables associated with the \objective\ constraints. This distinction will help define the dual feasibility constraints in the inverse formulations. 
\section{Proofs of Statements}

\begin{proof} {Proof of Proposition \ref{prop:IO_Feas}.}
Let $\ba_j$ be the j$^\text{th}$ row of $\bA$ and $\besmall_j$ be the j$^\text{th}$ unit vector (the set of \desirable\ constraints is non-empty) for some $j \in \cJ$. Then, $\exists \hat{\bx} \in \Omega$ for which the equality $\ba_j \hat{\bx} = b_j$ holds (given the assumption that no constraint is redundant). Defining $\bc_0=\ba'_j$,  $\by_0=\besmall _j$, $\bz_0 = \hat{\bx}$ and $\epsilon_0^k = \bx^k - \hat{\bx} \ \forall k \in \mathcal{K}$, the solution  $(\bc_0, \by_0, \bzero, \left \{ \epsilon_0^1, \cdots, \epsilon_0^K\right \}, \bz_0)$ is feasible for $\IO$ since \eqref{IOPrimalFeasblility1} and \eqref{IOPrimalFeasblility2} hold by the feasibility of $\hat{\bx}$, \eqref{IOStrongDual} holds by the assumption $\ba_j \hat{\bx} = b_j$ and the definition of $\bc_0$ and $\by_0$, and remaining constraints are also satisfied directly from the definition of the parameters.\Halmos 
\end{proof}

\begin{proof} {Proof of Proposition \ref{prop:SIO_regularization}.}
Assume to the contrary that for a feasible solution $({\bc}, {\by}, {\bu}, {\bE},{\bz})$ for $\IO$, we have ${\bc} = \bzero$ and therefore,  ${\by}' \bA = \bzero$. We have two cases:
\begin{itemize}
    \item \textbf{Case 1:} ${y}_j =0 \ \forall j \in \cJ$ which is a contradiction to the constraint \eqref{IORegularization}. 
    \item \textbf{Case 2:} There is a convex combination of the rows of matrix $\bA$ such that $\by' \bA = \bzero$. We show that this is a contradiction to the feasible set $\Omega$ being full dimensional. Considering that $\Omega$ does not have redundant constraints, let $\hat{\by}$ be the non-zero elements in ${\by}$ and let $\bD$ be the index set of these elements. Note that $2 \leq \left | \bD \right | \leq n$. Also, let $\hat{\bA}$ and $\hat{\bb}$ be sub-matrices of $\bA$ and $\bb$ containing rows indexed by $\bD$ and  let $\ba'_i$ be the $i^{th}$ row of $\hat{\bA}$ and $b_i$ be the $i^{th}$ element of right-hand-side vector $\hat{\bb}$ for $i \in \left\{ 1, \hdots, \left| \bD \right| \right\}$. We have the following:
    \begin{equation} \label{theorm1proof:eq1}
        \hat{y}_1 \ba'_1 + ... + \hat{y}_{\left| \bD \right|} \ba'_{\left| \bD \right|} = 0, \quad \hat{y}_i >0 \quad \forall i \in \left\{ 1, \hdots, \left| \bD \right| \right\} \\
         \Rightarrow \hat{y}_1 \ba'_1 + ... + \hat{y}_{{\left| \bD \right|}-1} \ba'_{{\left| \bD \right|}-1} = -\hat{y}_{\left| \bD \right|} \ba'_{\left| \bD \right|}. 
    \end{equation}
    Note that $\forall \bx \in \Omega$, we have the following:
    \begin{equation} 
        \left.\begin{matrix}
        \ba'_1 \bx \geq b_1\\ 
        \vdots \\ 
        \ba'_{{\left| \bD \right|}-1} \bx \geq b_{{\left| \bD \right|}-1}
        \end{matrix}\right\} \Rightarrow 
        (\hat{y}_1 \ba'_1 + ... + \hat{y}_{{\left| \bD \right|}-1} \ba'_{{\left| \bD \right|}-1}) \bx \geq \hat{y}_1 b_1 + ... + \hat{y}_{{\left| \bD \right|}-1} b_{{\left| \bD \right|}-1}
    \end{equation} 
    Which yields the following using \eqref{theorm1proof:eq1}:
    \begin{equation} \label{theorm1proof:eq2}
         -\hat{y}_{\left| \bD \right|} \ba'_{\left| \bD \right|} \bx \geq \sum_{i=1}^{{\left| \bD \right|}-1}\hat{y}_i b_i.
    \end{equation}     
    
    Additionally, $\forall \bx \in \Omega$:
    \begin{equation} \label{theorm1proof:eq3}
    \ba'_{{\left| \bD \right|}} \bx \geq b_{{\left| \bD \right|}} \Rightarrow -\hat{y}_{\left| \bD \right|} \ba'_{\left| \bD \right|} \bx \leq -\hat{y}_{\left| \bD \right|} b_{{\left| \bD \right|}}.    
    \end{equation}
    However, since strong duality conditions hold for solution $\bz$, from the complementary slackness conditions, we have the following:
    \begin{equation} \label{theorm1proof:eq4}
        \left.\begin{matrix}
        \ba'_1 \bz = b_1\\ 
        \vdots \\ 
        \ba'_{\left| \bD \right|} \bz = b_{\left| \bD \right|}
        \end{matrix}\right\} \Rightarrow 
        (\sum_{i=1}^{{\left| \bD \right|}} \hat{y}_{i} \ba'_{i}) \bz = \sum_{i=1}^{{\left| \bD \right|}}\hat{y}_i b_i.
    \end{equation} 
    However, we know that $\sum_{i=1}^{{\left| \bD \right|}} \hat{y}_{i} \ba'_{i}=0$. Therefore, $\sum_{i=1}^{{\left| \bD \right|}}\hat{y}_i b'_i = 0$ and we get the important equation $\sum_{i=1}^{{\left| \bD \right|}-1}\hat{y}_i b_i = -\hat{y}_{\left| \bD \right|} b_{\left| \bD \right|}$. Combining this equation with equality \eqref{theorm1proof:eq4} yields $-\hat{y}_{\left| \bD \right|} \ba'_{\left| \bD \right|} \bx \geq -\hat{y}_{\left| \bD \right|} b_{\left| \bD \right|}$. Considering this inequality with inequality \eqref{theorm1proof:eq3} yields $-\hat{y}_{\left| \bD \right|} \ba'_{\left| \bD \right|} \bx = -\hat{y}_{\left| \bD \right|} b_{{\left| \bD \right|}}$, which considering that $\hat{y}$ is non-zero, results in $\ba'_{\left| \bD \right|} \bx =  b_{{\left| \bD \right|}}$. This means that $\ba'_{\left| \bD \right|} \bx \geq  b_{{\left| \bD \right|}}$ is an implicit inequality for $\Omega$ which is a contradiction to $\Omega$ being full dimensional (for more information on implicit equalities refer to Theorem 3.17 from \cite{conforti2014integer}).
\Halmos
\end{itemize}
\end{proof}

\begin{proof} {Proof of Theorem \ref{Theorem_IL}.}
Since $(\bc^*, \by^*, \bar{\by}^*, \bE^*, \bz^*)$ is optimal for $\IO$, $\bz^*$ satisfies primal feasibility, dual feasibility, and strong duality conditions for $\FO (\bc^*, \Omega)$ and as such $\bz^*$ is optimal for $\FO (\bc^*, \Omega)$ and thus $\bz^* \in \Omega^{opt}(\bc^*)$. Since $\Omega^{opt}=\bigcup_{c\neq 0} \Omega^{opt}(\bc)$, $\bz^* \in \Omega^{opt}$.
\Halmos
\end{proof}

\begin{proof} {Proof of Proposition \ref{theorem1}.}
Since the objective functions of formulations \eqref{MFD} and the $\IO$ model are the same, it suffices to show that $(\ba'_j$, $e_j$,$\bzero$, $\bE^*$, $\bz^*)$ is optimal for $\IO$. It can be readily seen that $(\ba'_j$, $e_j$,$\bzero$, $\bE^*$, $\bz^*)$ is indeed feasible for $\IO$. Now assume to the contrary that there is a solution ($\bc$, $\by$,$\bzero$, $\bE$, $\bz$) for $\IO$ such that the objective value of $\IO$ for this solution is strictly less than $\cD_{min}$. However, from previous results, at least one \desirable\ constraint is binding at $\bz$. We denote this constraint as constraint $j$. As a result, we have $\ba_j  \bz = b_j$ which is exactly \eqref{ILjStrongDual}. \eqref{ILjOnePoint} is also easily seen to hold for $\bz$. Therefore, $\bz$ is feasible for formulation \eqref{MFD} which is a contradiction with $(\ba'_{j_{min}}$, $e_{j_{min}}$,$\bzero$, $\bE^*$, $\bz^*)$ being the solution with the minimum objective value among all solutions to different instances of formulations \eqref{MFD}. Therefore, $(\ba'_{j_{min}}$, $e_{j_{min}}$,$\bzero$, $\bE^*$, $\bz^*)$ is optimal for $\IO$.\Halmos
\end{proof}


\begin{proof} {Proof of Proposition \ref{MGIL_feas}}
We have the following arguments for each part:
\begin{enumerate}[(a)]
    \item Let $\ba_j$ be the j$^\text{th}$ row of $\bA$ and $\besmall_j$ be the j$^\text{th}$ unit vector (the set of acceptable constraints is non-empty). There exists an FO feasible point $\hat{\bx}$ for which the equality $\ba_j \hat{\bx} = b_j$ holds (given the assumption that no constraint is redundant). Defining $\bv_0=\besmall_j$, $\bz_0 = \hat{\bx}$ and $\epsilon_0^k = \bx^k - \hat{\bx} \ \forall k \in \mathcal{K}$, the solution  $(\bv_0, \left \{ \epsilon_0^1, \cdots, \epsilon_0^K\right \}, \bz_0)$ is feasible for $\MGIL(\bX,\Omega,r,\bP)$ with $r=1$ since: \eqref{GILStrongDual} and \eqref{GILPrimal Feasblility2} hold by the feasibility of $\hat{\bx}$ and the assumption $\ba_j \hat{\bx} = b_j$ and the definition of $\bv_0$, and remaining constraints are also satisfied directly from the definition of the parameters.
    
    \item Let $\bx^0$ be an extreme point of the feasible region of $\FO$ which binds exactly $n$ \desirable\ constraints. Let $B \subseteq \cJ$ be the set of indices of \desirable\ constraints that are binding at $\bx^0$. For $r_0 \in \left \{ 1,...,n \right \}$, let $B_{r_0} \subseteq B$ be a set of indices of $r_0$ binding constraints at $\bx^0$ and let $\bv_{r_0}$ be a vector such that $v_j = 1$ $\forall j \in B_{r_0}$ and $v_j = 0$ otherwise for $j \in \cJ$. Also, let $\bE_0$ be the matrix of distances from $\bx^0$ to all the observations. Then,  ($\bv_{r_0}$, $\bE_0$, $\bx^0$) is feasible for $\MGIL$ with $r = r_0$.
    
    \item If $\Omega$ has at least one face $\cF$ with $r_0$ binding \desirable\ constraints, then $\exists \bx^0$ that binds $r_0$ binding \desirable\ constraints. The arguments that $\MGIL$ is feasible for $r \in \left \{ 1,...,r_0 \right \}$ are similar to part 2.
    
    \item Let $\bz^*_{r_1}$ and $\bz^*_{r_2}$ be the optimal solutions of $\MGIL(\bX,\Omega,r_1,\bP)$ and $\MGIL(\bX,\Omega,r_2,\bP)$ respectively. We show that $\bz^*_{r_2}$ is a feasible solution for $\MGIL(\bX,\Omega,r_1,\bP)$. This can be easily seen by choosing $r_1$ number of binding \desirable\ constraints of $\bz^*_{r_2}$. Let $B$ be a subset of indices of binding constraints at $\bz^*_{r_2}$ with $\left | B \right | = r_1$ and let $\bv$ be a vector such that $v_j = 1$ $\forall j \in B$ and $v_j = 0$ otherwise for $j \in \cJ$. Also, let $\bE_{r_2}$ be the matrix of the distances of $\bz^*_{r_2}$ to the observations. It can be seen that ($\bv$, $\bE_{r_2}$, $\bz^*_{r_2}$) is feasible for $\MGIL(\bX,\Omega,r_1,\bP)$. Therefore, $\cD^*_{r_1} \leq \cD^*_{r_2}$.

\end{enumerate}
\end{proof}

\begin{proof} {Proof of Theorem \ref{Theorem2}.}
Let $\bc \in cone (\ba_t : t \in T)$. We have $\bc = \sum_{t \in T} \lambda_t \ba_t$ with $\lambda_t \geq 0$. It suffices to show that the dual feasibility and strong duality conditions hold for $\bz^*$ in $\FO(\bc,\Omega)$.  We have $\bc' \bz^* = (\sum_{t \in T} \lambda_t \ba_t') \bz^* = \sum_{t \in T}\lambda_t b_t$ with $\lambda_t \geq 0$ $\forall t \in T$. Let $\by$ be such that $y_j = \lambda_t$ $\forall t \in T$ and $y_j = 0$ otherwise. We have $\bA' \by = \bc$ and $\by \geq 0$ and $\bc' \bz^* = \bb' \by$.\Halmos
\end{proof}

\begin{proof} {Proof of Theorem \ref{rem:GILlinktoIL}}
$\Rightarrow$) Let $(\bc_{\IO}$, $\by_{\IO}$, $\bE_{\IO}$, $\bz_{\IO})$ optimal for $\IO$. We construct a feasible solution for $\MGIL$ based on the optimal solution of $\IO$. For $j \in \cJ$, let $v_j$ be the $j$th element of $\mathbf{v}_{\MGIL}$ such that $v_j = 1$ if $y_j > 0$ where $y_j$ is the $j$th element of $\by_{\IO}$ and zero otherwise. Then, it is easy to verify that $(\mathbf{v}_{\MGIL}$ $\bE_{\IO}$, $\bz_{\IO})$ is feasible for $\MGIL(\bX,\Omega,r = 1, \bP = \emptyset)$. If $(\mathbf{v}_{\MGIL}$ $\bE_{\IO}$, $\bz_{\IO})$ is optimal for $\MGIL(\bX,\Omega,r = 1, \bP = \emptyset)$, we are done. So assume to the contrary that there exists some solution  $(\hat{\mathbf{v}}$ $\hat{\bE}$, $\hat{\bz})$ optimal for $\MGIL(\bX,\Omega,r = 1, \bP = \emptyset)$ such that $\cD (\hat{\bE}) < \cD (\bE_{\IO})$. However, by similar arguments, another feasible solution to $\IO$ can be constructed based on this optimal solution to $\MGIL(\bX,\Omega,r = 1, \bP = \emptyset)$. But this is a contradiction to $(\bc_{\IO}$, $\by_{\IO}$, $\bE_{\IO}$, $\bz_{\IO})$ being optimal for $\IO$.

$\Leftarrow$) Showing that it is possible to construct feasible solutions from a solution of $\MGIL(\bX,\Omega,r = 1, \bP = \emptyset)$ to $\IO$ and vice versa in the previous part, this part can also be proven following similar arguments. \Halmos
\end{proof}

\begin{proof} {Proof of Theorem \ref{Thm:GGIL_Optimal}.}
For the sequence of solutions $\bz_1 , \hdots , \bz_L$, we have $\cJ_{\bz_1} \subseteq \hdots \subseteq \cJ_{\bz_L}$. As such, $\forall i \in \left \{ 1,\hdots,L \right \}$, $\bz_i$ binds the \desirable\ constraints indexed by $\cJ_{\bz_1}$. Therefore, $\forall i \in \left \{ 1,\hdots,L \right \}$, $\bz_i \in \mathcal{F} = \left \{ \bz \in \Omega | \ \ba_j \bz = b_j \ \forall j \in \cJ_{\bz_1} \right \}$
The results follow from Proposition \ref{prop:GGIL_Face}. Note that for each $\bz_i$ where $ i \in \left \{ 1,\hdots,L \right \}$, $\cJ_{z_1} \subseteq \hdots \subseteq \cJ_{z_i}$ and therefore, due to Theorem \ref{Theorem2}, the results hold. \Halmos
\end{proof}

\section{Diet Recommendation Problem Additional Data and Results} \label{Appendix:diet_data}
This section includes more information on the input data for the diet recommendation problem and complete figures indicating $\IO$ solutions for all food types and nutrients. The intake data include more than 5,000 different food types. Given the large number of food types, we bundled them into 38 broad food groups for ease of interpretation and to make the learned diets more tractable. This categorization is done based on the food codes from USDA. Table \ref{Table:food_groups} shows the grouping developed for the dataset and the average serving size of each food item in grams. Table \ref{Table:dash_diet_recommendations_servings} illustrates the recommendations of the DASH diet in terms of the number of servings of each food group for different diets with distinct calorie targets. Since the DASH diet recommendations are in servings, Table \ref{Table:dash_diet_recommendations_servings} provides additional details about a typical sample of each food group along with the corresponding amount in one serving size.  We utilize the food samples from Table \ref{Table:dash_diet_recommendations_servings}, the nutritional data from USDA, and the recommended amounts from the DASH eating plan to calculate the required bounds on nutrients. These bounds can serve as the right-hand side vector for constraints in linear optimization settings. Table \ref{Table:SubsetNutrients} depicts these amounts for the 1,600 and 2,000 calories target diets.
Figures \ref{Fig:ILvsMILallfoods} and \ref{Fig:ILvsMILnuts_all} are extended versions of Figures \ref{Fig:ILvsMILfoods} and \ref{Fig:ILvsMILnuts} showing the results of implementing inverse learning models for all food types and all nutrients. One additional detail to note is that in this extended figure of nutrients, the lower bound of sugar is binding for the $r=1$ of $\MGIL$ and in increasing values of $r$, it becomes non-binding. This is due to the fact that not only the binding constraints chosen by the binary variables are forced to remain binding and any other constraint that becomes binding in addition to the ones chosen by the binary variables might become non-binding in subsequent runs with higher $r$ values. With the increasing applications and importance of optimization and learning models, the ease of access to reliable, accurate, and interpretable datasets has become paramount. We provided the details of a large-scale open-access dataset on dietary behaviors and attributes of individuals in this section and have provided access to the data at the following address https://github.com/CSSEHealthcare/InverseLearning. We hope that the presence of such a dataset can help the researchers in different data-driven approaches to evaluate proposed methods and get meaningful insights.

\iftrue
\begin{figure}[t]
\begin{center}
\includegraphics[width =0.9 \linewidth]{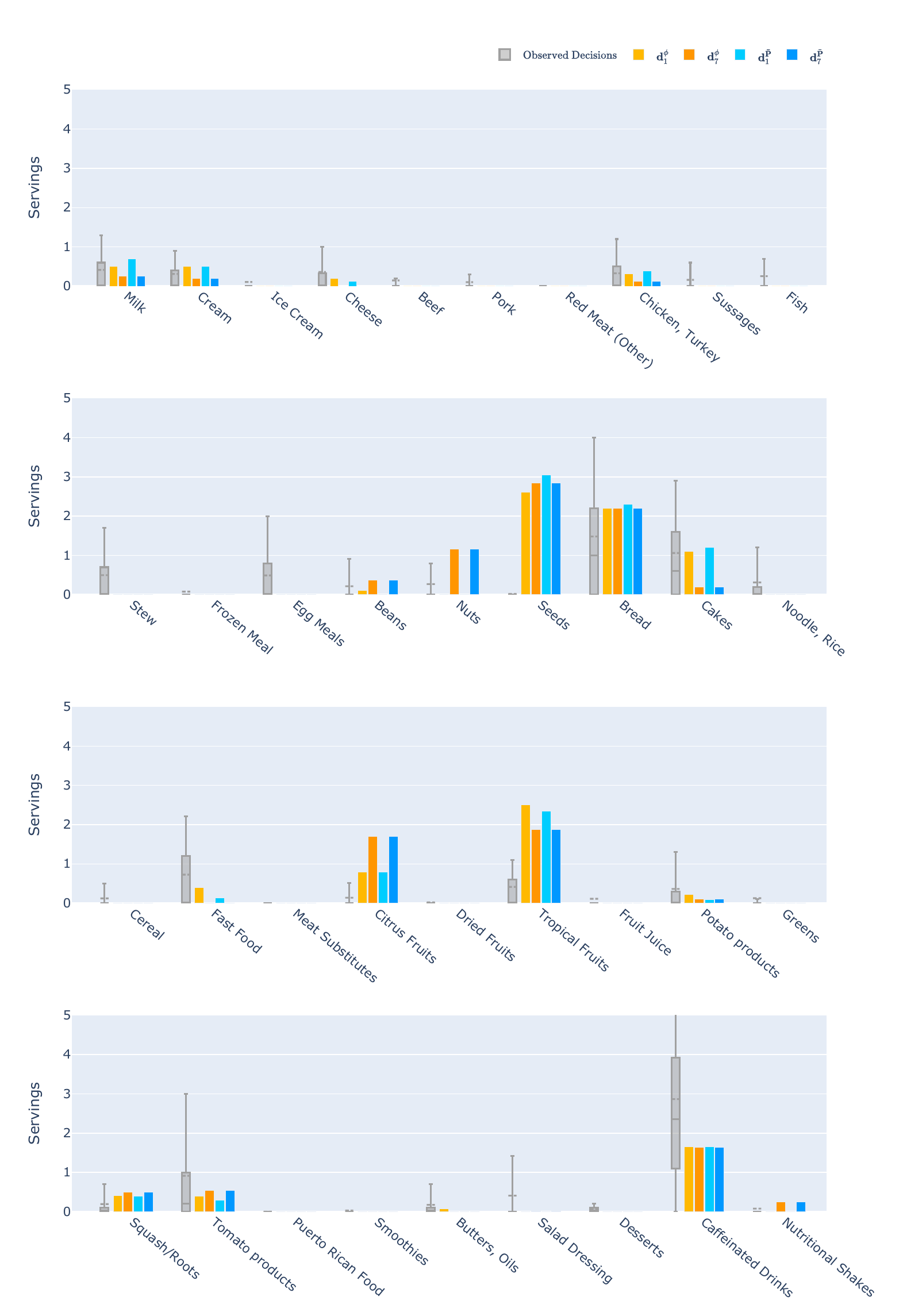}
\caption{\footnotesize Comparison of recommended diets by $\MGIL$ with different values for $r$ and a set of 460 observations for all food types. } \label{Fig:ILvsMILallfoods}
\end{center}
\end{figure}
\fi

\iftrue
\begin{figure}[t]
\begin{center}
\includegraphics[width =1 \linewidth]{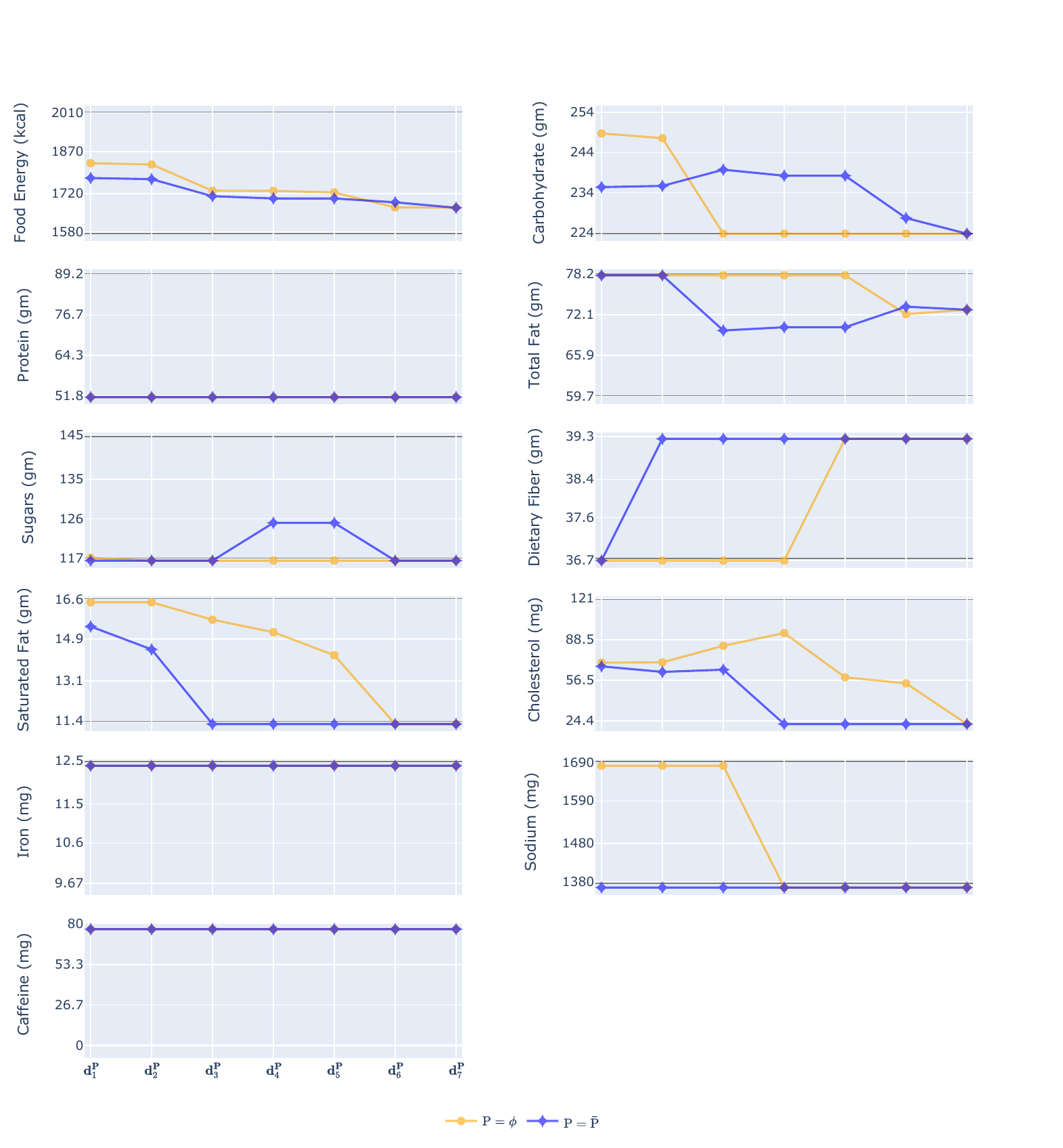}
\caption{\footnotesize Comparison of nutrients of recommended diets by $\GIL$ and $\MGIL$ for different values of $r$ for all nutrients in the model. (Triv.: \objective\ constraint, Rel.: \desirable\ constraint, Pref.: \pref\ constraint)} \label{Fig:ILvsMILnuts_all}
\end{center}
\end{figure}
\fi

\begin{table}[] 
\small
        \caption{Food groups and their respective serving sizes in grams}
        \label{Table:food_groups}
\begin{tabular}{>{}p{0.24\textwidth}|>{}p{0.55\textwidth}|>{\centering\arraybackslash}p{0.15\textwidth}}
Group Name              & Description                                                                                            & Serving Size (g) \\ \hline \hline
Milk                    & milk, soy milk, almond milk, chocolate milk, yogurt, baby food, infant   formula                         & 244              \\
Cream                   & Cream, sour cream                                                                                      & 32               \\
Ice Cream               & all types of ice cream                                                                                 & 130              \\
Cheese                  & all types of cheese                                                                                    & 32               \\
Beef                    & ground beef, steaks (cooked, boiled, grilled or raw)                                                   & 65               \\
Pork                    & chops of pork, cured pork, bacon (cooked, boiled, grilled or raw)                                      & 84               \\
Red Meat (Other)        & lamb, goat, veal, venison  (cooked, boiled, grilled or raw)                                            & 85               \\
Chicken, Turkey          & all types of chicken, turkey, duck (cooked, boiled, grilled or raw)                                   & 110              \\
Sausages                & beef or red meat by-products, bologna, sausages, salami, ham    (cooked, boiled, grilled or raw)      & 100              \\
Fish                    & all types of fish,                                                                                     & 85               \\
Stew                    & stew meals containing meat (or substitutes), rice, vegetables                                          & 140              \\
Frozen Meals            & frozen meal (containing meat and vegetables)                                                           & 312              \\
Egg Meals               & egg meals, egg omelets and substitutes                                                                 & 50               \\
Beans                   & all types of beans (cooked, boiled, baked, raw)                                                                                     & 130           \\
Nuts                    & all types of nuts                                                                                      & 28.35            \\
Seeds                   & all types of seeds                                                                                     & 30               \\
Bread                   & all types of bread                                                                                     & 25               \\
Cakes, Biscuits, Pancakes & cakes, cookies, pies, pancakes, waffles                                                                & 56               \\
Noodle, Rice             & macaroni, noodle, pasta, rice                                                                          & 176              \\
Cereal                  & all types of cereals                                                                                   & 55               \\
Fast Foods              & burrito, taco, enchilada, pizza, lasagna                                                               & 198              \\
Meat Substitutes        & meat substitute that are  cereal-   or vegetable protein-based                                         & 100              \\
Citrus Fruits           & grapefruits, lemons, oranges                                                                           & 236              \\
Dried Fruits            & all types of dried fruit                                                                               & 28.3            \\
Tropical Fruits         & apples, apricots, avocados, bananas, cantaloupes, cherries, figs, grapes,   mangoes, pears, pineapples & 182              \\
Fruit Juice             & All types of fruit juice                                                                               & 249              \\
Potato products                & potatoes (fried, cooked)                                                                                           & 117              \\
Greens                  & beet greens, collards, cress, romaine, greens, spinach                                                 & 38               \\
Squash/Roots            & carrots, pumpkins, squash, sweet potatoes                                                               & 72               \\
Tomato products                 & tomato, salsa containing tomatoes, tomato byproducts                                                   & 123              \\
Vegetables              & raw vegetables                                                                                         & 120              \\
Puerto Rican Food      & Puerto Rican style food                                                                                & 250              \\
Smoothies               & fruit and vegetable smoothies                                                                          & 233              \\
Butter, Oils             & butters, oils                                                                                          & 14.2             \\
Salad Dressing          & all types of salad dressing                                                                            & 14               \\
Desserts                & sugars, desserts, toppings                                                                             & 200              \\
Caffeinated Drinks      & coffees, soda drinks, iced teas                                                                        & 240              \\
Nutritional Shakes                  & Nutritional shakes, Energy Drinks,   Protein Powders                                                                                  & 166   \\    
\hline
\end{tabular}
\end{table}

\begin{table}[]
\small
        \caption{Food categories and their recommended number of servings for different targets based on the DASH diet \citep{dash_diet_2020}}
        \label{Table:dash_diet_recommendations_servings}
\begin{tabular}{>{}p{0.3\textwidth}|>{}p{0.1\textwidth}|>{}p{0.08\textwidth}|>{}p{0.08\textwidth}|>{}p{0.08\textwidth}|>{}p{0.08\textwidth}|>{}p{0.075\textwidth}|>{\arraybackslash}p{0.075\textwidth}} 
&\multicolumn{7}{c}{Diet Target}\\\cline{2-8}
Food Category                          & 1,200    \ \  Calories               & 1,400 Calories              & 1,600 Calories               & 1,800 Calories               & 2,000 Calories               & 2,600 Calories         & 3,100 Calories\\
\hline \hline
Grains                             & 4–5                & 5–6                & 6                  & 6                  & 6–8                & 10–11        & 12–13        \\
Vegetables                          & 3–4                & 3–4                & 3–4                & 4–5                & 4–5                & 5–6          & 6            \\
Fruits                              & 3–4                & 4                  & 4                  & 4–5                & 4–5                & 5–6          & 6            \\
Fat-free or low-fat dairy products & 2–3                & 2–3                & 2–3                & 2–3                & 2–3                & 3            & 3–4          \\
Lean meats, poultry, and fish       & $\leq$ 3          & $\leq$3–4        & $\leq$3–4         & $\leq$ 6         & $\leq$6          & $\leq$6     & 6–9          \\
Nuts, seeds, and legumes            & 3/week         & 3/week         & 3–4/week       & 4/week         & 4–5/week       & 1            & 1            \\
Fats and oils                      & 1                  & 1                  & 2                  & 2–3                & 2–3                & 3            & 4            \\
Sweets and added sugars             & $\leq$ 3/week & $\leq$3/week & $\leq$3/week & $\leq$5/week & $\leq$5/week & $\leq$2           & $\leq$2           \\
Maximum sodium limit(mg/day)               & 2,300        & 2,300        & 2,300        & 2,300       & 2,300        & 2,300  & 2,300 
\end{tabular}
\end{table}

\begin{table}[]
\small
        \caption{Food categories and their respective serving sizes in grams}
        \label{Table:dash_diet_servings}
\begin{tabular}{ll}
Food Category                            & Serving Size (Example)                                                        \\
\hline \hline
Grains                              & 1 slice of whole-grain bread                                                 \\
Vegetables                            & 1 cup (about 30 grams) of raw, leafy green vegetables like   spinach or kale \\
Fruits                                & 1 medium apple                                                               \\
Fat-free   or low-fat dairy products & 1 cup (240 ml) of low-fat milk                                               \\
Lean   meats, poultry, and fish       & 1 ounce (28 grams) of cooked meat, chicken or fish                           \\
Nuts,   seeds, and legumes            & 1/3 cup (50 grams) of nuts                                                   \\
Fats   and oils                     & 1 teaspoon (5 ml) of vegetable oil                                           \\
Sweets   and added sugars             & 1 cup (240 ml) of lemonade                                                  
\end{tabular}
\end{table}

\end{document}